\documentclass{aastex631}

\graphicspath{{./}{figures/}}
\usepackage{amsmath}
\usepackage{bbm}
\usepackage{mathrsfs}
\usepackage{mathtools}
\usepackage{amsthm}
\newtheorem{theorem}{Theorem}
\newtheorem{lemma}{Lemma}
\theoremstyle{remark}

\usepackage{booktabs}
\usepackage{multirow}
\setcitestyle{aysep={,}}

\allowdisplaybreaks

\begin{document}

\title{Asymptotic properties of resampling-based processes for the average treatment effect in observational studies with competing risks}

\correspondingauthor{Jasmin Rühl}
\email{jasmin.ruehl@math.uni-augsburg.de}

\author{Jasmin Rühl}
\affiliation{Department of Mathematical Statistics and Artificial Intelligence in Medicine \\
University of Augsburg \\
Augsburg, Germany}

\author{Sarah Friedrich}
\affiliation{Department of Mathematical Statistics and Artificial Intelligence in Medicine \\
University of Augsburg \\
Augsburg, Germany}
\affiliation{Centre for Advanced Analytics and Predictive Sciences (CAAPS) \\
University of Augsburg \\
Augsburg, Germany}


\begin{abstract}
In observational studies with time-to-event outcomes, the g-formula can be used to estimate a treatment effect in the presence of confounding factors. However, the asymptotic distribution of the corresponding stochastic process is complicated and thus not suitable for deriving confidence intervals or time-simultaneous confidence bands for the average treatment effect. A common remedy are resampling-based approximations, with Efron's nonparametric bootstrap being the standard tool in practice. We investigate the large sample properties of three different resampling approaches and prove their asymptotic validity in a setting with time-to-event data subject to competing risks.
\end{abstract}


\section{Introduction}

In observational studies, comparisons between treatment groups are complicated by the potentially unequal distribution of confounding factors. A prominent idea to tackle this issue is the potential outcomes approach, which models the mean outcome in a hypothetical world where all study participants are subject to the same intervention \citep{rubin1974estimating, hernan2020causal}. This article focuses on the comparison between two treatment groups, where the outcome is a right-censored time-to-event that is possibly subject to competing risks. Our parameter of interest is the average treatment effect (ATE) at time $t$, defined as the difference between the absolute risks of the event of interest in the exposure groups. Following \citet{ozenne2020on}, we model cause-specific hazards by Cox regression models to estimate the absolute risk of the event of interest \citep{benichou1990estimates, ozenne2017riskRegression}. \\
For thorough statistical inference, confidence intervals and time-simultaneous confidence bands for the ATE provide additional insight. Inference is usually based on Efron's nonparametric bootstrap, since the stochastic processes involved in the estimation are rather complex, though \citep{efron1981censored, stensrud2020separable, ryalen2020causal, neumann2016covariate}. \citet{ozenne2020on} present an alternative approach, which is based on the influence function and the resampling scheme developed by \citet{scheike2008flexible}. Another popular resampling technique for time-to-event data is a martingale-based wild bootstrap. This method has successfully been applied in different situations, e.g. in non-causal investigations that cover Cox proportional hazards models \citep{lin1993checking, lin1994confidence} or competing risks \citep{lin1997non-parametric}. Extensions that improve the performance for small sample sizes have been discussed by \citet{beyersmann2013weak}, \citet{dobler2014bootstrapping} as well as \citet{dobler2017non}. What is more, the wild bootstrap has been shown to perform superior to the classical bootstrap in several situations, in particular when the data involve dependencies \citep{niessl2023statistical, ruehl2022general}.

In this paper, we derive a martingale representation of the stochastic process characterizing the asymptotic behavior of the ATE. Based on this representation, we provide proofs of the asymptotic validity of three resampling approaches: Efron's bootstrap, a resampling  method based on the influence function and the wild bootstrap. Thus, the main contribution of this paper is to fill the gap between theory and practice and provide the missing proofs that justify the application of resampling techniques in the situation discussed here. \\
The remainder of this manuscript is organized as follows: 
Section~\ref{sec:setting} introduces the notation, the competing risks setting and the parameter of interest. In Section~\ref{sec:asymptotics}, we investigate the asymptotic behavior of the estimated average treatment effect and derive a martingale-based representation of the corresponding stochastic process. Section~\ref{sec:resampling} contains the proofs of the asymptotic validity of the three resampling approaches. We close with a discussion in Section~\ref{sec:discussion}.

\section{Setting and Notation \label{sec:setting}}

Consider an independent and identically distributed (i.i.d.) data sample of the form $\smash{\{(T_i \!\land\! C_i, D_i, A_i, \boldsymbol{Z_i})\}_{i \in \{1, \dots, n\}}}$. The first vector element denotes the time of an individual's first event ($T$) or censoring ($C$), whichever occurs earlier. In a setting with $K$ competing causes of failure, the indicator $D$ may assume values in ${\{1, \dots, K\}}$ according to the type of event that is observed, whereas for censored observations, ${D = 0}$. The data moreover include the treatment indicator ${A \in \{0,1\}}$ as well as the covariate vector $\smash{\boldsymbol{Z} \in \mathbb{R}^p}$. We suppose that there are no ties and that ${T \perp \!\!\! \perp C \mid (A, \boldsymbol{Z})}$. Besides,the covariate values in $\boldsymbol{Z}$ should be bounded. 

Let without loss of generality ${D = 1}$ refer to the event of interest and consider the 
potential cumulative incidence function $\smash{F_1^a(t) = P(T^a \leq t, \, D^a = 1)}$ under treatment $a$. We characterize the average treatment effect by the relation $\smash{ATE(t) = \mathbb{E}\left(F_1^1(t) - F_1^0(t)\right)}$, with time $t$ ranging between 0 and $\tau$, the terminal time of the study. If the identifiability conditions of exchangeability, positivity and consistency apply \citep[cf.][section~I.3]{hernan2020causal}, the g-formula suggests
\[\widehat{ATE}(t) = \frac{1}{n} \sum_{i=1}^n \left(\hat{F}_1(t \mid A = 1, \boldsymbol{Z_i}) - \hat{F}_1(t \mid A = 0, \boldsymbol{Z_i})\right)\]
as an estimator of the average treatment effect \citep{ozenne2020on}. One possible way to obtain $\smash{\hat{F}_1(t \mid A, \boldsymbol{Z})}$ involves fitting cause-specific Cox models with covariates $A$ and $\smash{\boldsymbol{Z^{(k)}} \in \mathbb{R}^{p_k}}$ for each event type ${k \in \{1, \dots, K\}}$. This yields estimated cumulative hazards of the following form:
\begin{gather*}
\hat{\Lambda}_k(t \mid a, \boldsymbol{z^{(k)}}) = \hat{\Lambda}_{0k}(t) \exp(\hat{\beta}_{kA} a + \skew{3.5}\hat{\boldsymbol{\beta}}_{\boldsymbol{kZ}}^T \boldsymbol{z^{(k)}}), \\
\text{with } \hat{\Lambda}_{0k}(t) = \int_0^t \frac{dN_k(s)}{\sum_{i=1}^n Y_i(s) \exp(\hat{\beta}_{kA} A_i + \skew{3.5}\hat{\boldsymbol{\beta}}_{\boldsymbol{kZ}}^T \boldsymbol{Z_i^{(k)}})}.
\end{gather*}
(We use $Z$ instead of $\smash{Z^{(k)}}$ hereafter, as the cause specificity of the covariates follows from the context.) The vector $\smash{\skew{3.5}\hat{\boldsymbol{\beta}}_{\boldsymbol{k}} = (\hat{\beta}_{kA}, \skew{3.5}\hat{\boldsymbol{\beta}}_{\boldsymbol{kZ}}^T)^T \in \mathbb{R}^{p_k+1}}$ results from the Cox regression and combines the estimated coefficients for treatment and covariates. Apart from that, $\smash{\hat{\Lambda}_{0k}}(t)$ is the Breslow estimator of the cumulative baseline hazard \citep{breslow1972contribution}, which depends on the counting process $\smash{N_k(t) = \sum_{i=1}^n N_{ki}(t)}$, with $\smash{N_{ki}(t) = \mathbbm{1}\{T_i \!\land\! C_i \leq t, \, D_i = k\}}$, and the at-risk indicator $\smash{Y_i(t) = \mathbbm{1}\{T_i \!\land\! C_i \geq t\}}$, ${i \in \{1, \dots, n\}}$. Note that we use ${\mathbbm{1}\{\cdot\}}$ to denote the indicator function. The estimator of the cumulative incidence finally results by plugging $\smash{\Lambda_k(t \mid a, \boldsymbol{z})}$, into the formula
\[\hat{F}_1(t \mid a, \boldsymbol{z}) = \int_0^t \hat{S}(t \mid a, \boldsymbol{z}) \, d\hat{\Lambda}_1(s \mid a, \boldsymbol{z}),\]
where $\smash{\hat{S}(t \mid a, \boldsymbol{z}) = \exp\left(-\sum_{k=1}^K \hat{\Lambda}_k(s \mid a, \boldsymbol{z})\right)}$ approximates the survival probability ${P(T > t \mid a, \boldsymbol{z})}$ for a given treatment $a$ and covariate vector $\boldsymbol{z}$.
 
\section{Asymptotic distribution of the estimated average treatment effect \label{sec:asymptotics}}

In order to investigate the asymptotic behavior of $\smash{\widehat{ATE}}$, we study the process $\smash{U_n(t) = \sqrt{n} \, (\widehat{ATE}(t) - ATE(t))}$ and its properties as ${n \to \infty}$. Arguments similar to those used by \citet{cheng1998prediction} show that the limiting distribution of $U_n$ may be represented in terms of martingales. This is an important finding that facilitates further inferences on the large-sample properties of the process. Before we commence with the proof, it is necessary to define several functions and variables, however. \\
Consider
\[S^{(r)}(\boldsymbol{\beta_k}, t) = \frac{1}{n} \sum_{i=1}^n Y_i(t) \exp(\beta_{kA} A_i + \boldsymbol{\beta}_{\boldsymbol{kZ}}^T \boldsymbol{Z_i}) \left((A_i, \boldsymbol{Z}_{\boldsymbol{i}}^T)^T\right)^{\otimes r}\!\!, \quad r \in \{0,1,2\},\]
(where $\smash{\boldsymbol{v}^{\otimes 0} = 1}$, $\smash{\boldsymbol{v}^{\otimes 1} = \boldsymbol{v}}$, $\smash{\boldsymbol{v}^{\otimes 2} = \boldsymbol{v} \boldsymbol{v}^T}$ for a column vector $\boldsymbol{v}$), and the corresponding expectations $\smash{s^{(r)}(\boldsymbol{\beta_k}, t) = \mathbb{E}\left(S^{(r)}(\boldsymbol{\beta_k}, t)\right)}$, ${r \in \{0, 1, 2\}}$, as well as
\[\boldsymbol{E}(\boldsymbol{\beta_k}, t) = \frac{\boldsymbol{S^{(1)}}(\boldsymbol{\beta_k}, t)}{S^{(0)}(\boldsymbol{\beta_k}, t)}\] 
with $\smash{\boldsymbol{e}(\boldsymbol{\beta_k}, t) = \frac{\boldsymbol{s^{(1)}}(\boldsymbol{\beta_k}, t)}{s^{(0)}(\boldsymbol{\beta_k}, t)}}$ for ${k \in \{1, \dots, K\}}$. The positive definite matrix $\smash{\mathbf{\Sigma_k}}$ is given by
\[\mathbf{\Sigma_k} = \int_0^{\tau} \left(\frac{\mathbf{s^{(2)}}(\boldsymbol{\beta_{0k}}, u)}{s^{(0)}(\boldsymbol{\beta_{0k}}, u)} - \left(\boldsymbol{e}(\boldsymbol{\beta_{0k}}, u)\right)^{\otimes 2}\right) s^{(0)}(\boldsymbol{\beta_{0k}}, u) \, d\Lambda_{0k}(u),\]
supposing that the Cox model applies with true vector of regression coefficients $\smash{\boldsymbol{\beta_{0k}}}$ for cause $k$. Let further
\begin{gather*}
\boldsymbol{h_k}(t \mid a, \boldsymbol{z}) = \int_0^t \left((a, \boldsymbol{z}^T)^T - \boldsymbol{e}(\boldsymbol{\beta_{0k}}, u)\right) \, d\Lambda_{k}(u \mid a, \boldsymbol{z}), \\[-0.2cm]
\boldsymbol{\varphi_1}(t \mid a, \boldsymbol{z}) = \int_0^t S(u- \mid a, \boldsymbol{z}) \, d\boldsymbol{h_1}(u \mid a, \boldsymbol{z}), \\[-0.2cm]
\boldsymbol{\psi_{1k}}(t \mid a, \boldsymbol{z}) = \int_0^t \left(F_1(t \mid a, \boldsymbol{z}) - F_1(u \mid a, \boldsymbol{z})\right) \, d\boldsymbol{h_k}(u \mid a, \boldsymbol{z}).
\end{gather*}
Eventually, we define the functions 
\[H_{k1i}(u,t) = \frac{\tilde{H}_{k1}(u, t)}{\sqrt{n} \!\cdot\! S^{(0)}(\boldsymbol{\beta_{0k}}, u)} \quad \text{and} \quad H_{k2i}(u,t) = \frac{1}{\sqrt{n}} \, \left(\boldsymbol{\tilde{H}_{k2}}(t)\right)^T \! \mathbf{\Sigma_k}^{\!\!\!\!-1} \left((A_i, \boldsymbol{Z}_{\boldsymbol{i}}^T)^T - \boldsymbol{E}(\boldsymbol{\beta_{0k}}, u)\right),\]
${k \in \{1, \dots, K\}}$, ${i \in \{1, \dots, n\}}$, ${u \leq t}$, with 
\begin{gather*}
\begin{align*}
\tilde{H}_{11}(u,t) &= \frac{1}{n} \sum_{i=1}^n \left(\left(S(u- \mid A=1, \boldsymbol{Z_i}) - F_1(t \mid A=1, \boldsymbol{Z_i}) + F_1(u \mid A=1, \boldsymbol{Z_i})\right) \exp(\beta_{01A} + \boldsymbol{\beta}_{\boldsymbol{01Z}}^T \boldsymbol{Z_i}) \right. \\[-0.3cm]
&\hphantom{= \frac{1}{n} \sum_{i=1}^n \left( i \right.} \left. - \left(S(u- \mid A=0, \boldsymbol{Z_i}) - F_1(t \mid A=0, \boldsymbol{Z_i}) + F_1(u \mid A=0, \boldsymbol{Z_i})\right) \exp(\boldsymbol{\beta}_{\boldsymbol{01Z}}^T \boldsymbol{Z_i})\right),
\end{align*} \\
\begin{align*}
\tilde{H}_{k1}(u,t) &= \frac{1}{n} \sum_{i=1}^n \left(\left(F_1(t \mid A=0, \boldsymbol{Z_i}) - F_1(u \mid A=0, \boldsymbol{Z_i})\right) \exp(\boldsymbol{\beta}_{\boldsymbol{0kZ}}^T \boldsymbol{Z_i}) \right. \\[-0.3cm]
&\hphantom{= \frac{1}{n} \sum_{i=1}^n \left( i \right.} \left. - \left(F_1(t \mid A=1, \boldsymbol{Z_i}) - F_1(u \mid A=1, \boldsymbol{Z_i})\right) \exp(\beta_{0kA} + \boldsymbol{\beta}_{\boldsymbol{0kZ}}^T \boldsymbol{Z_i}) \right)
\end{align*} \hfill \quad (k \in \{2, \dots, K\})
\end{gather*}
and
\begin{gather*}
\boldsymbol{\tilde{H}_{12}}(t) = \frac{1}{n} \sum_{i=1}^n \left(\left(\boldsymbol{\varphi_1}(t \mid A=1, \boldsymbol{Z_i}) - \boldsymbol{\psi_{11}}(t \mid A=1, \boldsymbol{Z_i})\right) - \left(\boldsymbol{\varphi_1}(t \mid A=0, \boldsymbol{Z_i}) - \boldsymbol{\psi_{11}}(t \mid A=0, \boldsymbol{Z_i})\right)\right), \\
\boldsymbol{\tilde{H}_{k2}}(t) = \frac{1}{n} \sum_{i=1}^n \left(\boldsymbol{\psi_{1k}}(t \mid A=0, \boldsymbol{Z_i}) - \boldsymbol{\psi_{1k}}(t \mid A=1, \boldsymbol{Z_i})\right) \tag*{$(k \in \{2, \dots, K\})$.}
\end{gather*}

\begin{lemma} \label{lemma_martingale_representation}
For the process
\[\tilde{U}_n(t) = \sum_{k=1}^K \sum_{i=1}^n \left(\int_0^t H_{k1i}(s,t) \, dM_{ki}(s) + \int_0^\tau H_{k2i}(s,t) \, dM_{ki}(s)\right),\]
with $\smash{M_{ki}(t) = N_{ki}(t) - \int_0^t Y_i(s) \, d\Lambda_k(s \mid A_i, \boldsymbol{Z_i})}$, ${k \in \{1, \dots, K\}}$, ${i \in \{1, \dots, n\}}$, it holds that 
\[U_n(t) = \tilde{U}_n(t) + o_p(1).\]
\end{lemma}
Note that $\smash{M_{ki}}$ is a martingale relative to the history $\smash{\left(\mathscr{F}_t\right)_{t \geq 0}}$ that is generated by the data observed until a given time,
i.e.~$\smash{\mathbb{E}\left(dM_{ki}(t) \mid \mathscr{F}_{t-}\right) = 0}$. 
\begin{proof}[Proof]
By the strong law of large numbers, we have
\begin{align*}
U_n(t) &= \frac{\sqrt{n}}{n} \sum_{i=1}^n \left(\int_0^t \hat{S}(u- \mid A=1, \boldsymbol{Z_i}) \exp(\hat{\beta}_{1A} + \skew{3.5}\hat{\boldsymbol{\beta}}_{\boldsymbol{1Z}}^T \boldsymbol{Z_i}) \, d\hat{\Lambda}_{01}(u) - \int_0^t S(u- \mid A=1, \boldsymbol{Z_i}) \exp(\beta_{01A} + \boldsymbol{\beta}_{\boldsymbol{01Z}}^T \boldsymbol{Z_i}) \, d\Lambda_{01}(u)\right) \\[-0.1cm]
&\hphantom{= i} - \frac{\sqrt{n}}{n} \sum_{i=1}^n \left(\int_0^t \hat{S}(u- \mid A=0, \boldsymbol{Z_i}) \exp(\skew{3.5}\hat{\boldsymbol{\beta}}_{\boldsymbol{1Z}}^T \boldsymbol{Z_i}) \, d\hat{\Lambda}_{01}(u) - \int_0^t S(u- \mid A=0, \boldsymbol{Z_i}) \exp(\boldsymbol{\beta}_{\boldsymbol{01Z}}^T \boldsymbol{Z_i}) \, d\Lambda_{01}(u)\right) + o_p(1) \\
&= \frac{\sqrt{n}}{n} \sum_{i=1}^n \left(\int_0^t \left(\hat{S}(u- \mid A=1, \boldsymbol{Z_i}) - S(u- \mid A=1, \boldsymbol{Z_i})\right) \exp(\hat{\beta}_{1A} + \skew{3.5}\hat{\boldsymbol{\beta}}_{\boldsymbol{1Z}}^T \boldsymbol{Z_i}) \, d \hat{\Lambda}_{01}(u) \right. \\[-0.2cm]
&\hphantom{= \frac{\sqrt{n}}{n} \sum_{i=1}^n \left( i \right.} \left. + \int_0^t S(u- \mid A=1, \boldsymbol{Z_i}) \,d\left(\exp(\hat{\beta}_{1A} + \skew{3.5}\hat{\boldsymbol{\beta}}_{\boldsymbol{1Z}}^T \boldsymbol{Z_i}) \hat{\Lambda}_{01}(u) - \exp(\beta_{01A} + \boldsymbol{\beta}_{\boldsymbol{01Z}}^T \boldsymbol{Z_i}) \Lambda_{01}(u)\right)\right) \\[-0.1cm]
&\hphantom{= i} - \frac{\sqrt{n}}{n} \sum_{i=1}^n \left(\int_0^t \left(\hat{S}(u- \mid A=0, \boldsymbol{Z_i}) - S(u- \mid A=0, \boldsymbol{Z_i})\right) \exp(\skew{3.5}\hat{\boldsymbol{\beta}}_{\boldsymbol{1Z}}^T \boldsymbol{Z_i}) \, d \hat{\Lambda}_{01}(u) \right. \\[-0.2cm]
&\hphantom{= i -\frac{\sqrt{n}}{n} \sum_{i=1}^n \left( i \right.} \left. + \int_0^t S(u- \mid A=0, \boldsymbol{Z_i}) \,d\left(\exp(\skew{3.5}\hat{\boldsymbol{\beta}}_{\boldsymbol{1Z}}^T \boldsymbol{Z_i}) \hat{\Lambda}_{01}(u) - \exp(\boldsymbol{\beta}_{\boldsymbol{01Z}}^T \boldsymbol{Z_i}) \Lambda_{01}(u)\right)\right) + o_p(1).
\end{align*}
\citet{lin1994confidence} showed that $\smash{\sqrt{n} \left(\hat{\Lambda}_k(t \mid a, \boldsymbol{z}) - \Lambda_k(t \mid a, \boldsymbol{z})\right) = \tilde{W}_k(t \mid a, \boldsymbol{z}) + o_p(1)}$ for the martingale expression 
\[\tilde{W}_k(t \mid a, \boldsymbol{z}) = \frac{1}{\sqrt{n}} \sum_{i=1}^n \left(\int_0^t \frac{\exp(\beta_{0kA} \!\cdot\! a + \boldsymbol{\beta}_{\boldsymbol{0kZ}}^T \boldsymbol{z})}{S^{(0)}(\boldsymbol{\beta_{0k}}, u)} \, dM_{ki}(u) + \left(\boldsymbol{h_k}(t \mid a, \boldsymbol{z})\right)^T  \mathbf{\Sigma_k}^{\!\!\!\!-1} \! \int_0^{\tau} \left((A_i, \boldsymbol{Z}_{\boldsymbol{i}}^T)^T - \boldsymbol{E}(\boldsymbol{\beta_{0k}}, u)\right) \, dM_{ki}(u)\right).\]
Thus, exploiting the (uniform) consistency of $\smash{\skew{3.5}\hat{\boldsymbol{\beta}}_{\boldsymbol{1}}}$ and $\smash{\hat{\Lambda}_{01}}$ (\citealp{tsiatis1981a}; \citealp[pp.~361--362]{kosorok2008introduction}) and using a first-order Taylor approximation of ${f: x \mapsto \exp(-x)}$ around $\smash{x = \sum_{k=1}^K \exp(\beta_{0kA} \!\cdot\! a + \boldsymbol{\beta}_{\boldsymbol{0kZ}}^T \boldsymbol{z}) \Lambda_{0k}(t)}$ (which yields $\smash{\hat{S}(t- \mid a, \boldsymbol{z}) - S(t- \mid a, \boldsymbol{z}) = -\frac{1}{\sqrt{n}} S(t- \mid a, \boldsymbol{z}) \sum_{k=1}^K \tilde{W}_k(t \mid a, \boldsymbol{z}) + o_p(1)}$), we find that
\begin{align*}
U_n(t) &= \frac{1}{n} \sum_{i=1}^n \left(\int_0^t S(u- \mid A=1, \boldsymbol{Z_i}) \, d\tilde{W}_1(u \mid A=1, \boldsymbol{Z_i}) - \sum_{k=1}^K \int_0^t \tilde{W}_k(u \mid A=1, \boldsymbol{Z_i}) \, dF_1(u \mid A=1, \boldsymbol{Z_i})\right) \\[-0.1cm]
&\hphantom{= i} - \frac{1}{n} \sum_{i=1}^n \left(\int_0^t S(u- \mid A=0, \boldsymbol{Z_i}) \, d\tilde{W}_1(u \mid A=0, \boldsymbol{Z_i}) - \sum_{k=1}^K \int_0^t \tilde{W}_k(u \mid A=0, \boldsymbol{Z_i}) \, dF_1(u \mid A=0, \boldsymbol{Z_i})\right) + o_p(1) \\
&= \frac{1}{n} \sum_{i=1}^n \left(\int_0^t S(u- \mid A=1, \boldsymbol{Z_i}) \, d\tilde{W}_1(u \mid A=1, \boldsymbol{Z_i}) \right. \\[-0.1cm]
&\hphantom{= \frac{1}{n} \sum_{i=1}^n \left( i \right.} \left. - \sum_{k=1}^K \int_0^t \left(F_1(t \mid A=1, \boldsymbol{Z_i}) - F_1(u \mid A=1, \boldsymbol{Z_i})\right) \, d\tilde{W}_k(u \mid A=1, \boldsymbol{Z_i}) \right. \\[-0.1cm]
&\hphantom{= \frac{1}{n} \sum_{i=1}^n \left( i \right.} \left. - \left(\int_0^t S(u- \mid A=0, \boldsymbol{Z_i}) \, d\tilde{W}_1(u \mid A=0, \boldsymbol{Z_i}) \right. \right. \\[-0.1cm]
&\hphantom{= \frac{1}{n} \sum_{i=1}^n \left(i -\left( i \right. \right.} \left. \left. - \sum_{k=1}^K \int_0^t \left(F_1(t \mid A=0, \boldsymbol{Z_i}) - F_1(u \mid A=0, \boldsymbol{Z_i})\right) \, d\tilde{W}_k(u \mid A=1, \boldsymbol{Z_i})\right) \right) + o_p(1).
\end{align*}
The last equivalence follows from integration by parts, since 
\[\int_0^t \tilde{W}_k(u \mid a, \boldsymbol{z}) \, dF_1(u \mid a, \boldsymbol{z}) = \tilde{W}_k(t \mid a, \boldsymbol{z}) F_1(t \mid a, \boldsymbol{z}) - \int_0^t F_1(u \mid a, \boldsymbol{z}) \, d\tilde{W}_k(u \mid a, \boldsymbol{z}).\]
Finally, by inserting the definition of $\smash{\tilde{W}_k}$ and reordering the terms, the result follows.
\end{proof}

The subsequent theorem characterizes the limiting distribution of $\smash{U_n}$ for fixed covariate vectors $\smash{\boldsymbol{Z_i}}$, ${i \in \{1, \dots, n\}}$:
\begin{theorem} \label{theorem:limit_U}
The process $\smash{U_n}$ converges weakly to a zero-mean Gaussian process (with covariance function $\smash{\xi(t_1, t_2) = \sum_{k=1}^K \xi^{(k)}(t_1, t_2)}$,
\[\xi^{(k)}(t_1, t_2) = \int_0^{t_1 \land t_2} \tilde{H}_{k1}(u, t_1) \tilde{H}_{k1}(u, t_2) \frac{d \Lambda_{0k}(u)}{s^{(0)}(\boldsymbol{\beta_{0k}}, u)} + \left(\boldsymbol{\tilde{H}_{k2}}(t_1)\right)^T \! \mathbf{\Sigma_k}^{\!\!\!\!-1} \, \boldsymbol{\tilde{H}_{k2}}(t_2),\]
on the Skorokhod space ${\mathcal{D}[0, \tau]}$.
\end{theorem}
\begin{proof}[Proof]
Lemma~\ref{lemma_martingale_representation} implies that it is sufficient to consider the limiting distribution of $\smash{\tilde{U}_n}$. \\
For distinct causes ${k \neq l}$, the counting processes $\smash{N_{ki}}$ and $\smash{N_{li}}$ cannot jump both, which is why the martingales $\smash{M_{ki}(t)}$ and $\smash{M_{li}(t)}$ are orthogonal. Moreover, ${\forall k \in \{1, \dots, K\}}$, 
\begin{align*}
<\sum_{i=1}^n \int_0^\cdot & \frac{1}{\sqrt{n} \!\cdot\! S^{(0)}(\boldsymbol{\beta_{0k}}, u)} \, dM_{ki}(u), \  \sum_{i=1}^n \int_0^\cdot \frac{1}{\sqrt{n}} \left((A_i, \boldsymbol{Z}_{\boldsymbol{i}}^T)^T - \boldsymbol{E}(\boldsymbol{\beta_{0k}}, u)\right) \,dM_{ki}(u) >(t) \\ 
= & \; \frac{1}{n} \sum_{i=1}^n \int_0^t \frac{1}{S^{(0)}(\boldsymbol{\beta_{0k}}, u)} \left((A_i, \boldsymbol{Z}_{\boldsymbol{i}}^T)^T - \boldsymbol{E}(\boldsymbol{\beta_{0k}}, u)\right) Y_i(u) \exp(\beta_{0kA} \!\cdot\! A_i + \boldsymbol{\beta}_{\boldsymbol{0kZ}}^T \boldsymbol{Z_i}) \, d \Lambda_{0k}(u) \\
= & \; \int_0^t \frac{1}{S^{(0)}(\boldsymbol{\beta_{0k}}, u)} \left(\boldsymbol{S^{(1)}}(\boldsymbol{\beta_{0k}}, u) - \boldsymbol{E}(\boldsymbol{\beta_{0k}}, u) S^{(0)}(\boldsymbol{\beta_{0k}}, u)\right) \, d \Lambda_{0k}(u) = \boldsymbol{0}.
\end{align*}
This means that $\smash{\sum_{i=1}^n \int_0^t H_{k1i}(u,t) \, dM_{ki}(u)}$ and $\smash{\sum_{i=1}^n \int_0^t H_{k2i}(u, t) \, dM_{ki}(u)}$ are orthogonal as well. \citet[pp.~498--501]{andersen1993statistical} furthermore showed that 
\[\frac{1}{\sqrt{n}} \sum_{i=1}^n \int_0^\tau \left((A_i, \boldsymbol{Z}_{\boldsymbol{i}}^T)^T - \boldsymbol{E}(\boldsymbol{\beta_{0k}}, u)\right) \, dM_{ki}(u) \  \overset{\mathscr{D}}{\longrightarrow} \  \mathcal{N}(\boldsymbol{0}, \mathbf{\Sigma_k})\]
as $n$ tends to infinity (where $\mathcal{N}$ symbolizes the normal distribution), and since $\smash{\boldsymbol{\varphi_1}}$ and $\smash{\boldsymbol{\psi_{1k}}}$ are deterministic functions for given covariates $\smash{\boldsymbol{Z_i}}$, the second summand of $\smash{\tilde{U}_n}$ is likewise asymptotically normal with mean zero. \\
It therefore only remains to consider the first summand. Note that ${\forall k \in \{1, \dots, K\}}$, the processes $\smash{\tilde{H}_{k1}(u,t)}$ are deterministic, continuous in $u$ and bounded for fixed covariates $\smash{\boldsymbol{Z_i}}$. In particular, 
\begin{gather*}
\left|\smash{\tilde{H}_{k1}}(u,t)\right| \leq \left(\exp(\beta_{0kA}) + 1\right) \max_{1 \leq i \leq n} \exp(\boldsymbol{\beta}_{\boldsymbol{0kZ}}^T \boldsymbol{Z_i}), \\
\left|\left(\smash{\tilde{H}_{k1}}(u,t)\right)^2\right| \leq \left(\exp(2 \beta_{0kA}) + 2 \exp(\beta_{0kA}) + 1\right) \max_{1 \leq i,j \leq n} \exp(\boldsymbol{\beta}_{\boldsymbol{0kZ}}^T (\boldsymbol{Z_i} + \boldsymbol{Z_j}))
\end{gather*}
for ${u \leq t}$. The strong law of large numbers further suggests that $\smash{S^{(0)}(\boldsymbol{\beta_k}, t)}$ converges to $\smash{s^{(0)}(\boldsymbol{\beta_k}, t)}$ almost surely for any ${t \in [0,\tau]}$, $\smash{\boldsymbol{\beta_k} \in \mathbb{R}^{p_k+1}}$. If we suppose that ${P(Y_i(\tau) > 0) > 0}$ ${\forall i \in \{1, \dots, n\}}$ (or also some less stringent constraints, see \citealp[section~8.4]{fleming2005counting}), this convergence is uniform on $\smash{\mathcal{B}_k \times [0, \tau]}$, where $\smash{\mathcal{B}_k}$ is a neighborhood of $\smash{\boldsymbol{\beta_{0k}}}$. Besides, $\smash{s^{(0)}}$ is bounded away from zero on $\smash{\mathcal{B}_k \times [0, \tau]}$. The conditions of Rebolledo's martingale central limit theorem \citep[Theorem~II.5.1]{andersen1993statistical} are thus fulfilled, and we conclude that $\smash{\tilde{U}_n}$ converges weakly to a zero-mean Gaussian process on ${\mathcal{D}[0,\tau]}$.

For the covariance function, one finds that
\begin{align*}
\tilde{\xi}(t_1, t_2) &= \sum_{k=1}^K \left( \sum_{i=1}^n \int_0^{t_1 \land t_2} \frac{\tilde{H}_{k1}(u, t_1) \tilde{H}_{k1}(u, t_2)}{n \left(\smash{S^{(0)}(\boldsymbol{\beta_{0k}},u)}\right)^2} Y_i(u) \exp(\beta_{0kA} \!\cdot\! A_i + \boldsymbol{\beta}_{\boldsymbol{0kZ}}^T \boldsymbol{Z_i}) \, d \Lambda_{0k}(u) \right. \\[-0.1cm]
&\hphantom{= \sum_{k=1}^K \left( i \right.} + \sum_{i=1}^n \int_0^\tau \left. \frac{1}{n} \left(\boldsymbol{\tilde{H}_{k2}}(t_1)\right)^T \! \mathbf{\Sigma_k}^{\!\!\!\!-1} \left((A_i, \boldsymbol{Z}_{\boldsymbol{i}}^T)^T - \boldsymbol{E}(\boldsymbol{\beta_{0k}}, u)\right) \left(\boldsymbol{\tilde{H}_{k2}}(t_2)\right)^T \! \mathbf{\Sigma_k}^{\!\!\!\!-1} \left((A_i, \boldsymbol{Z}_{\boldsymbol{i}}^T)^T - \boldsymbol{E}(\boldsymbol{\beta_{0k}}, u)\right) \right. \\[-0.3cm]
&\hphantom{= \sum_{k=1}^K \left( i + \sum_{i=1}^n \int_0^\tau i \right.} \left. \cdot Y_i(u) \exp(\beta_{0kA} \!\cdot\! A_i + \boldsymbol{\beta}_{\boldsymbol{0kZ}}^T \boldsymbol{Z_i}) \, d \Lambda_{0k}(u) \vphantom{\sum_{i=1}^n} \right) \\
&= \sum_{k=1}^K \left(\int_0^{t_1 \land t_2} \frac{\tilde{H}_{k1}(u, t_1) \tilde{H}_{k1}(u, t_2)}{S^{(0)}(\boldsymbol{\beta_{0k}},u)} \, d \Lambda_{0k}(u) \right. \\[-0.1cm]
&\hphantom{= \sum_{k=1}^K \left( i \right.} \left. + \left(\boldsymbol{\tilde{H}_{k2}}(t_1)\right)^T \! \mathbf{\Sigma_k}^{\!\!\!\!-1} \left(\int_0^\tau \frac{1}{n} \sum_{i=1}^n \left((A_i, \boldsymbol{Z}_{\boldsymbol{i}}^T)^T - \boldsymbol{E}(\boldsymbol{\beta_{0k}}, u)\right) \left((A_i, \boldsymbol{Z}_{\boldsymbol{i}}^T)^T - \boldsymbol{E}(\boldsymbol{\beta_{0k}}, u)\right)^T \right. \right. \\[-0.3cm]
&\hphantom{= \sum_{k=1}^K \left( i + \left(\boldsymbol{\tilde{H}_{k2}}(t_1)\right)^T \! \mathbf{\Sigma_k}^{\!\!\!\!-1} \left(\int_0^\tau \frac{1}{n} \sum_{i=1}^n i \right. \right.} \left. \left. \cdot Y_i(u) \exp(\beta_{0kA} \!\cdot\! A_i + \boldsymbol{\beta}_{\boldsymbol{0kZ}}^T \boldsymbol{Z_i}) \, d \Lambda_{0k}(u) \vphantom{\sum_{i=1}^n}\right) \left(\mathbf{\Sigma_k}^{\!\!\!\!-1}\right)^T \! \boldsymbol{\tilde{H}_{k2}}(t_2) \vphantom{\int_0^{t_1}}\right) \\
&= \sum_{k=1}^K \left(\int_0^{t_1 \land t_2} \frac{\tilde{H}_{k1}(u, t_1) \tilde{H}_{k1}(u, t_2)}{S^{(0)}(\boldsymbol{\beta_{0k}},u)} \, d \Lambda_{0k}(u) \right. \\[-0.1cm]
&\hphantom{= \sum_{k=1}^K \left( i \right.} \left. + \left(\boldsymbol{\tilde{H}_{k2}}(t_1)\right)^T \! \mathbf{\Sigma_k}^{\!\!\!\!-1} \left(\int_0^\tau \left(\mathbf{S^{(2)}}(\boldsymbol{\beta_{0k}},u) - \boldsymbol{S^{(1)}}(\boldsymbol{\beta_{0k}},u) \left(\boldsymbol{E}(\boldsymbol{\beta_{0k}},u)\right)^T\right) \, d \Lambda_{0k}(u)\right) \mathbf{\Sigma_k}^{\!\!\!\!-1} \, \boldsymbol{\tilde{H}_{k2}}(t_2)\right) \\
&\underset{n \to \infty}{\longrightarrow} \sum_{k=1}^K \left(\int_0^{t_1 \land t_2} \frac{\tilde{H}_{k1}(u, t_1) \tilde{H}_{k1}(u, t_2)}{s^{(0)}(\boldsymbol{\beta_{0k}},u)} \, d \Lambda_{0k}(u) + \left(\boldsymbol{\tilde{H}_{k2}}(t_1)\right)^T \! \mathbf{\Sigma_k}^{\!\!\!\!-1} \mathbf{\Sigma_k} \mathbf{\Sigma_k}^{\!\!\!\!-1} \, \boldsymbol{\tilde{H}_{k2}}(t_2) \right) = \xi(t_1, t_2),
\end{align*}
where the convergence in the last step follows by the strong law of large numbers and the continuous mapping theorem.
\end{proof}

\section{Resampling-based approximations \label{sec:resampling}}

The asymptotic distribution of $\smash{U_n(t)}$ is too complex to derive in practice, which is why resampling approaches are often used as a remedy to draw inferences on ${ATE(t)}$. In the following, we show the validity of three different methods.

\subsection{Efron's bootstrap}

Usually, confidence intervals and bands for the average treatment effect are constructed using the classical nonparametric bootstrap \citep{efron1981censored}. The main idea is to draw $n$ times with replacement from the data at hand and compute the desired statistical functional in the resulting bootstrap sample. This step is repeated multiple times, yielding a set of bootstrap estimators that provides information on the distribution of the underlying functional. Although this approach generally provides asymptotically valid outcomes, there are certain situations where it breaks down \citep{singh1981on, friedrich2017permuting}. To the best of our knowledge, a proof of the validity in the specific setting considered here is still pending.

\begin{theorem}
$\smash{U_n^*(t) = \sqrt{n} \left(\smash{\widehat{ATE}}^*(t) - \widehat{ATE}(t)\right)}$, with $\smash{\smash{\widehat{ATE}}^*(t)}$ being the estimated average treatment effect in the bootstrap sample, converges to the same limiting process as $\smash{U_n(t)}$ for almost all data samples $\smash{\{T_i \!\land\! C_i, D_i, A_i, \boldsymbol{Z_i}\}_{i \in \{1, \dots, n\}}}$ if $\smash{\inf_{u \in [0, \tau]} Y(u) \overset{P}{\to} \infty}$.
\end{theorem}
The superscript `$*$' is used here and in the following to indicate bootstrapped quantities.
\begin{proof}[Proof (Outline)]
Suppose that the given data were obtained on the probability space ${(\Omega, \mathcal{A}, P)}$. We first note that the general martingale arguments apply conditionally on ${\omega \in \Omega}$ for almost all $\omega$ \citep[cf.][]{akritas1986bootstrapping}. Let $\smash{\tau^* = \max_{1 \leq i \leq n}\{(T \!\land\! C)_i^*\}}$. The estimators $\smash{\hat{\Lambda}_{0k}}$ and $\smash{\hat{\Lambda}_k}$ calculated in the original data sample are the true (discontinuous) cumulative baseline hazard and cumulative hazard in the bootstrap sample, respectively. Moreover, $\smash{s^{(r), *}(\boldsymbol{\beta_k}, u) = S^{(r)}(\boldsymbol{\beta_k}, u)}$ as well as $\smash{\boldsymbol{e^{*}}(\boldsymbol{\beta_k}, u) = \boldsymbol{E}(\boldsymbol{\beta_k}, u)}$, which is easy to see if the bootstrap sample is represented with multinomial weights assigned to the original sample, e.g.
\[S^{(0), *}(\boldsymbol{\beta_k}, u) = \frac{1}{n} \sum_{i=1}^n w_i \, Y_i(u) \exp(\beta_{kA} A_i + \boldsymbol{\beta}_{\boldsymbol{kZ}}^T \boldsymbol{Z_i})\]
for $\smash{\boldsymbol{w} \sim \mathcal{M}ult(n, (1/n, \dots, 1/n)^T)}$. Thus, $\smash{\mathbf{\Sigma_k^*} = \hat{\mathbf{\Sigma}}_\mathbf{k}}$, $\smash{\boldsymbol{h_k^*} = \boldsymbol{\hat{h}_k}}$, $\smash{\boldsymbol{\varphi_1^*} = \boldsymbol{\hat{\varphi}_1^*}}$, $\smash{\boldsymbol{\psi_{1k}^*} = \boldsymbol{\hat{\psi}_{1k}}}$ and
\[M_{ki}^*(t) = w_{i} \left(N_{ki}(t) - \int_0^t Y_i(u) \, d\hat{\Lambda}_k(u \mid A_i, \boldsymbol{Z_i})\right),\]
${i \in \{1, \dots, n\}}$. Note that a discrete-time setting is considered here! We can now infer that $\smash{\boldsymbol{\hat{\beta}_k^*} \overset{P}{\to} \boldsymbol{\hat{\beta}_k}}$ and $\smash{\hat{\Lambda}_{0k}^* \overset{a.s.}{\to} \hat{\Lambda}_{0k}}$ on $\smash{[0, \tau^*]}$ as ${n \to \infty}$ by the considerations of \citet{prentice2003mixed}. Also,
\begin{align*}
\sqrt{n} \left(\hat{\Lambda}_k^*(t \mid a, \boldsymbol{z}) - \hat{\Lambda}_k(t \mid a, \boldsymbol{z})\right) &= \frac{1}{\sqrt{n}} \int_0^t \frac{\exp(\hat{\beta}_{kA} \!\cdot\! a + \skew{3.5}\hat{\boldsymbol{\beta}}_{\boldsymbol{kZ}}^T \boldsymbol{z})}{S^{(0), *}(\boldsymbol{\hat{\beta}_k}, u)} M_k^*(du) \\
&\hphantom{= i} + \frac{1}{\sqrt{n}} \left(\boldsymbol{\hat{h}_k}(t \mid a, \boldsymbol{z})\right)^T \! \hat{\mathbf{\Sigma}}_\mathbf{k}^{-1} \left(\sum_{i=1}^n \int_0^{\tau^*} \left(\left(A_i, \boldsymbol{Z}_{\boldsymbol{i}}^T\right)^T - \boldsymbol{E^*}(\boldsymbol{\hat{\beta}_k}, u)\right) M_{ki}^*(du)\right) + o_p(1),
\end{align*}
which can be concluded by the reasoning of \citet[proof of Theorem~VII.2.3]{andersen1993statistical}. These results provide the basis for proceeding in the same way as we did in the proof of Lemma~\ref{lemma_martingale_representation}. It follows that $\smash{U_n^*(t) = \tilde{U}_n^*(t) + o_p(1)}$,
with 
\[\tilde{U}_n^*(t) = \sum_{k=1}^K \sum_{i=1}^n \left(\int_0^t H_{k1i}^*(u,t) M_{ki}^*(du) + \int_0^{\tau^*} H_{k2i}^*(u,t) M_{ki}^*(du)\right),\] 
applying the definitions of $\smash{H_{k1i}}$ and $\smash{H_{k2i}}$ to the bootstrap sample.

Subsequently, we use similar arguments as in the proof of Theorem~\ref{theorem:limit_U}. One finds that
\begin{align*}
<\sum_{i=1}^n \int_0^\cdot & \frac{1}{\sqrt{n} \!\cdot\! S^{(0), *}(\boldsymbol{\hat{\beta}_k}, u)} M_{ki}^*(du), \sum_{i=1}^n \int_0^\cdot \frac{1}{\sqrt{n}} \left(\left(A_i, \boldsymbol{Z}_{\boldsymbol{i}}^T\right)^T - \boldsymbol{E^*}(\boldsymbol{\hat{\beta}_k}, u)\right) M_{ki}^*(du) >(t) \\
= & \; \frac{1}{n} \sum_{i=1}^m \int_0^t \frac{1}{S^{(0), *}(\boldsymbol{\hat{\beta}_k}, u)} \left(\left(A_i, \boldsymbol{Z}_{\boldsymbol{i}}^T\right)^T - \boldsymbol{E^*}(\boldsymbol{\hat{\beta}_k}, u)\right) w_i \, Y_i(u) \left(1 - \hat{\Lambda}_k(\Delta u \mid A_i, \boldsymbol{Z_i})\right) \hat{\Lambda}_k(du \mid A_i, \boldsymbol{Z_i}) \\
%
%
\longrightarrow & \; \int_0^t \frac{1}{S^{(0)}(\boldsymbol{\hat{\beta}_k}, u)} \left(\boldsymbol{E}(\boldsymbol{\hat{\beta}_k}, u) \, S^{(0)}_2(\boldsymbol{\hat{\beta}_k}, u) - \boldsymbol{S^{(1)}_2}(\boldsymbol{\hat{\beta}_k}, u)\right) \hat{\Lambda}_{0k} (\Delta u) \hat{\Lambda}_{0k}(du)
\end{align*}
with $\smash{S^{(r)}_2(\boldsymbol{\beta_k}, u) = \frac{1}{n} \sum_{i=1}^n Y_i(t) \exp(2 \beta_{kA} A_i + 2 \boldsymbol{\beta}_{\boldsymbol{kZ}}^T \boldsymbol{Z_i}) \left(\left(A_i, \boldsymbol{Z}_{\boldsymbol{i}}^T\right)^T\right)^{\otimes r}}$, ${r \in \{0,1,2\}}$. Also,
\begin{align*}
\frac{1}{\sqrt{n}} \sum_{i=1}^n \int_0^{\tau^*} \left(\left(A_i, \boldsymbol{Z}_{\boldsymbol{i}}^T\right)^T - \boldsymbol{E^*}(\boldsymbol{\hat{\beta}_k}, u)\right) M_{ki}^*(du) \overset{\mathscr{D}}{\longrightarrow} \mathcal{N}(\boldsymbol{0}, \hat{\mathbf{\Sigma}}_\mathbf{k} - \hat{\mathbf{\Sigma}}_\mathbf{k, 2})
\end{align*}
with
\begin{align*}
\hat{\mathbf{\Sigma}}_\mathbf{k, 2} &= \int_0^\tau \left(\mathbf{S^{(2)}_2}(\boldsymbol{\hat{\beta}_k}, u) - \boldsymbol{S^{(1)}_2}(\boldsymbol{\hat{\beta}_k}, u) \left(\boldsymbol{E}(\boldsymbol{\hat{\beta}_k}, u)\right)^T - \boldsymbol{E}(\boldsymbol{\hat{\beta}_k}, u) \left(\boldsymbol{S^{(1)}_2}(\boldsymbol{\hat{\beta}_k}, u)\right)^T  \right. \\[-0.1cm]
&\hphantom{= \int_0^\tau \left( i \right.} \left. +  \boldsymbol{E}(\boldsymbol{\hat{\beta}_k}, u) \left(\boldsymbol{E}(\boldsymbol{\hat{\beta}_k}, u)\right)^T \! S^{(0)}_2(\boldsymbol{\hat{\beta}_k}, u)\right) \hat{\Lambda}_{0k}(\Delta u) \hat{\Lambda}_{0k}(du)
\end{align*}
\citep[cf.][]{prentice2003mixed}, and lastly,
\begin{align*}
\sum_{i=1}^n \int_0^t \frac{\tilde{H}_{k1}^*(u, t)}{\sqrt{n} \!\cdot\! S^{(0), *}(\boldsymbol{\hat{\beta}_k}, u)} M_{ki}^*(du)
\end{align*}
converges weakly to a zero-mean Gaussian process with covariance function 
\[\xi_1^*(t_1, t_2) = \sum_{k=1}^K \int_0^{t_1 \land t_2} \frac{\tilde{H}_{k1}(u, t_1) \tilde{H}_{k1}(u, t_2)}{S^{(0)}(\boldsymbol{\hat{\beta}_k}, u)} \, \hat{\Lambda}_{0k}(du) - \sum_{k=1}^K \int_0^{t_1 \land t_2} \frac{\tilde{H}_{k1}(u, t_1) \tilde{H}_{k1}(u, t_2)}{\left(S^{(0)}(\boldsymbol{\hat{\beta}_k}, u)\right)^2} \, S^{(0)}_2(\boldsymbol{\hat{\beta}_k}, u) \, \hat{\Lambda}_{0k}(\Delta u) \hat{\Lambda}_{0k}(du)\]
as ${n \to \infty}$. \\
Since we assumed that there are no ties in the original sample,
\begin{align*}
S^{(r)}_2(\boldsymbol{\hat{\beta}_k}, u) \, \hat{\Lambda}_{0k}(\Delta u) \hat{\Lambda}_{0k}(du) &= \frac{1}{n} \sum_{i=1}^n Y_i(u) \exp(2 \hat{\beta}_{kA} A_i + 2 \skew{3.5}\hat{\boldsymbol{\beta}}_{\boldsymbol{kZ}}^T \boldsymbol{Z_i}) \left((A_i, \boldsymbol{Z}_{\boldsymbol{i}}^T)^T\right)^{\otimes r} \, \frac{\left(\Delta N_k(u)\right)^2}{\left(\sum_{i=1}^n Y_i(u) \exp(\hat{\beta}_{kA} A_i + \skew{3.5}\hat{\boldsymbol{\beta}}_{\boldsymbol{kZ}}^T \boldsymbol{Z_i})\right)^2} \\
&\leq \frac{ \max_{1 \leq i \leq n: \, Y_i(u) = 1}\{\exp(2 \hat{\beta}_{kA} A_i + 2 \skew{3.5}\hat{\boldsymbol{\beta}}_{\boldsymbol{kZ}}^T \boldsymbol{Z_i}) \left((A_i, \boldsymbol{Z}_{\boldsymbol{i}}^T)^T\right)^{\otimes r}\}}{Y(u)^2 \, \min_{1 \leq i \leq n: \, Y_i(u) = 1}\{\exp(2 \hat{\beta}_{kA} A_i + 2\skew{3.5}\hat{\boldsymbol{\beta}}_{\boldsymbol{kZ}}^T \boldsymbol{Z_i})\}}.
\end{align*}
Because of the boundedness of the covariates, all the terms involving $\smash{S^{(r)}_2(\boldsymbol{\hat{\beta}_k}, u) \hat{\Lambda}_{0k}(\Delta u) \hat{\Lambda}_{0k}(du)}$, ${r \in \{0, 1, 2\}}$, vanish as ${n \to \infty}$, and the proof is complete.
\end{proof}

\subsection{Influence function}

\citet{ozenne2020on} presented a second resampling technique based on the influence function of the average treatment effect. The idea proceeds from the functional delta method, which shows that
\[U_n(t) = \frac{1}{\sqrt{n}} \sum_{i=1}^n I\!F(t; \, T_i \!\land\! C_i, \! D_i, \! A_i, \! \boldsymbol{Z_i}) + o_P(1) \  \overset{\mathscr{D}}{\longrightarrow} \  \mathcal{N}\left(0, \int \left(I\!F(t; \, s, \! d, \! a, \! z)\right)^2 \, dP(s, \! d, \! a, \! z)\right).\]
To look up the definition of the influence function ${I\!F}$, refer to \citet{ozenne2020on, ozenne2017riskRegression}. The authors propose the resampling method described by \citet{scheike2008flexible} in order to approximate the distribution of the process while taking the dependence of the increments of $\smash{U_n}$ into account. This method is valid because of the asymptotic properties of $\smash{U_n}$ (see Theorem \ref{theorem:limit_U}). More specifically, one can imitate the limiting distribution of $\smash{U_n}$ by applying independent standard normal variables $\smash{G_1^{I\!F}, \dots, G_n^{I\!F}}$ and the plug-in estimator $\smash{\widehat{I\!F}}$ as follows:
\[\frac{1}{\sqrt{n}} \sum_{i=1}^n \widehat{I\!F}(t; T_i \!\land\! C_i, \! D_i, \! A_i, \! \boldsymbol{Z_i}) \cdot G_i^{I\!F}.\]
For further details, see also \citet[chapter~20]{vandervaart1998asymptotic}.

\subsection{Wild bootstrap}

With regard to the martingale representation of $\smash{U_n}$ from Lemma~\ref{lemma_martingale_representation}, a third approach results in accordance with the resampling scheme proposed by \citet{lin1993checking}. In short, one tries to emulate the distribution of the martingale increments $\smash{dM_{ki}}$, ${k \in \{1, \dots, K\}}$, ${i \in \{1, \dots, n\}}$ by generating random variates with asymptotically equal moments. The subsequent theorem sets out the conditions these variates need to fulfill in more detail. (Note the parallels to Theorem~1 in \citealp{dobler2017non}.)
\begin{theorem}
\label{theorem:U_hat}
Let $\smash{G_i^{W\!B}}$, ${i \in \{1, \dots, n\}}$, be random variables that satisfy the following conditions:
\begin{enumerate}
\item[(i)] $\smash{\sqrt{n} \, \max_{1 \leq i \leq n} \mathbb{E}\left(G_i^{W\!B} \mid \mathscr{F}_{\tau}\right) \overset{P}{\longrightarrow} 0}$;
\item[(ii)] $\smash{\max_{1 \leq i \leq n} \text{\normalfont{Var}}\left(G_i^{W\!B} \mid \mathscr{F}_{\tau}\right) \overset{P}{\longrightarrow} 1}$;
\item[(iii)] $\smash{\tfrac{1}{\sqrt{n}} \, \max_{1 \leq i \leq n} \mathbb{E}\left(\left(G_i^{W\!B}\right)^4 \mid \mathscr{F}_{\tau}\right) \overset{P}{\longrightarrow} 0}$;
\item[(iv)] $\smash{\mathcal{L}\left(G_i^{W\!B}, \, i \in \{1, \dots, n\} \mid \mathscr{F}_{\tau}\right) = \bigotimes_{i=1}^n \mathcal{L}\left(G_i^{W\!B} \mid \mathscr{F}_{\tau}\right)}$, \\
where $\mathcal{L}(\cdot \mid \mathscr{F}_{\tau})$ denotes the conditional distribution given $\mathscr{F}_{\tau}$ and $\otimes$ is the product measure;
\item[(v)] $\smash{\sum_{i=1}^n \mathbb{E}\left(\frac{\left(G_i^{W\!B} - \mathbb{E}\left(G_i^{W\!B} \mid \mathscr{F}_{\tau}\right)\right)^2}{\sum_{j=1}^n \left(\text{Var}\left(G_j^{W\!B} \mid \mathscr{F}_{\tau}\right)\right)^2} \cdot \mathbbm{1}\left\{\frac{\left(G_i^{W\!B} - \mathbb{E}\left(G_i^{W\!B} \mid \mathscr{F}_{\tau}\right)\right)^2}{\sum_{j=1}^n \left(\text{Var}\left(G_j^{W\!B} \mid \mathscr{F}_{\tau}\right)\right)^2} > \epsilon\right\} \mid \mathscr{F}_{\tau}\right) \overset{P}{\longrightarrow} 0 \  \forall \epsilon > 0}$.
\end{enumerate}
Then, conditional on the data,
\[\hat{U}_n(t) = \sum_{k=1}^K \sum_{i=1}^n \left(\hat{H}_{k1i}(T_i \!\land\! C_i, t) \, N_{ki}(t) G_i^{W\!B} + \hat{H}_{k2i}(T_i \!\land\! C_i, t) \, N_{ki}(\tau) G_i^{W\!B} \right)\]
converges weakly to the same process as $U_n$ on ${\mathcal{D}[0, \tau]}$.
\end{theorem}
The functions $\hat{H}_{k1i}$ and $\hat{H}_{k2i}$ are calculated by plugging appropriate sample estimates into the definitions of $H_{k1i}$ and $H_{k2i}$. It is easy to see that conditions (i) to (v) are fulfilled by independent standard normal multipliers $G_i^{W\!B}$, which corresponds to the original idea of \citet{lin1993checking}. Another option are independent centered unit Poisson multipliers, according to the suggestion of \citet{beyersmann2013weak}.

Before we can verify Theorem~\ref{theorem:U_hat}, several interim results are needed. The proofs of the following lemmas can be found in the appendix; the ideas are based on \citet{beyersmann2013weak} and \citet{dobler2017non}. \\
Consider the triangular arrays $\smash{\boldsymbol{X_{n,i}^{(k)}} = \left(\smash{X_{n,i}^{(k)}(t_1), \dots, X_{n,i}^{(k)}(t_l)}\right)^T}$, ${i \in \{1, \dots, n\}}$, ${k \in \{1, \dots, K\}}$, defined on the probability space $\smash{(\Omega_1, \mathcal{A}_1, P_1)}$, for $\smash{0 \leq t_1 \leq \dots \leq t_l \leq \tau}$, ${l \in \mathbb{N}}$, with
\[X_{n,i}^{(k)}(t) = \int_0^t \hat{\tilde{H}}_{k1}(u,t) \frac{dN_{ki}(u)}{\sqrt{n} \!\cdot\! S^{(0)}(\skew{3.5}\hat{\boldsymbol{\beta}}_{\boldsymbol{k}}, u)} + \int_0^\tau \frac{1}{\sqrt{n}} \, \left(\boldsymbol{\hat{\tilde{H}}_{k2}}(t)\right)^T \! \hat{\mathbf{\Sigma}}_\mathbf{k}^{-1} \left((A_i, \boldsymbol{Z}_{\boldsymbol{i}}^T)^T - \boldsymbol{E}(\hat{\boldsymbol{\beta}_k}, u)\right) \, dN_{ki}(u),\]
plug-in estimators $\smash{\hat{\tilde{H}}_{k1}}$, $\smash{\hat{\tilde{H}}_{k2}}$ and 
\[\hat{\mathbf{\Sigma}}_\mathbf{k} = \frac{1}{n} \sum_{i=1}^n \int_0^\tau \left(\frac{\mathbf{S^{(2)}}(\skew{3.5}\hat{\boldsymbol{\beta}}_{\boldsymbol{k}}, u)}{S^{(0)}(\skew{3.5}\hat{\boldsymbol{\beta}}_{\boldsymbol{k}}, u)} - \left(\frac{\boldsymbol{S^{(1)}}(\skew{3.5}\hat{\boldsymbol{\beta}}_{\boldsymbol{k}}, u)}{S^{(0)}(\skew{3.5}\hat{\boldsymbol{\beta}}_{\boldsymbol{k}}, u)}\right)^{\otimes 2} \right) \, dN_{ki}(u).\] 
Consequently, $\smash{\hat{U}_n(t) = \sum_{k=1}^K \sum_{i=1}^n G_i X_{n,i}^{(k)}(t)}$, with multipliers $\smash{G_i}$ defined on $\smash{(\Omega_2, \mathcal{A}_2, P_2)}$. (We generally consider the product probability space $\smash{(\Omega, \mathcal{A}, P) = (\Omega_1 \!\times\! \Omega_2, \mathcal{A}_1 \!\otimes\! \mathcal{A}_2, P_1 \!\times\! P_2)}$.)

\begin{lemma} \label{lemma_conditions_X}
The triangular arrays $\smash{\boldsymbol{X_{n,i}^{(k)}}}$ satisfy the following conditions for each ${k \in \{1, \dots, K\}}$: 
\begin{enumerate}
\item[(i)] $\smash{\max_{1 \leq i \leq n} || \boldsymbol{X_{n,i}^{(k)}} || \overset{P}{\longrightarrow} 0}$ (where ${||\cdot ||}$ denotes the Euclidian norm);
\item[(ii)] $\smash{\sum_{i=1}^n \boldsymbol{X_{n,i}^{(k)}} \left(\boldsymbol{X_{n,i}^{(k)}}\right)^T \overset{P}{\longrightarrow} \left(\xi^{(k)}(t_r, t_s)\right)_{1 \leq r, s \leq l}}$.
\end{enumerate}
\end{lemma}
\begin{lemma} \label{lemma_boundedness_X}
For time points $\smash{0 \leq t_r \leq t_s \leq \tau}$ and ${k \in \{1, \dots, K\}}$, 
\[\max_{1 \leq i \leq n} \left|X_{n,i}^{(k)}(t_s) - X_{n,i}^{(k)}(t_r)\right| \in O_p(n^{-1/2}),\]
where $\smash{O_p(a_n)}$ denotes asymptotic boundedness by $\smash{a_n}$ in probability.
\end{lemma}
Note that the bound in the Lemma~\ref{lemma_boundedness_X} is independent of the time points $\smash{t_r}$ and $\smash{t_s}$!
\begin{lemma} \label{lemma_boundedness_expectation_GX}
For time points $\smash{0 \leq t_q \leq t_r \leq t_s \leq \tau}$, causes ${k \in \{1, \dots, K\}}$ and the function
\begin{align*}
L_n^{(k)}(t) &= \frac{1}{n} \sum_{j=1}^n \left(\exp(2 \hat{\beta}_{kA}) + 2\exp(\hat{\beta}_{kA} + 1)\right) \exp(2 \skew{3.5}\hat{\boldsymbol{\beta}}_{\boldsymbol{kZ}}^T \boldsymbol{Z_j}) \int_0^t \frac{1}{n} \frac{dN_k(u)}{\left(S^{(0)}(\skew{3.5}\hat{\boldsymbol{\beta}}_{\boldsymbol{k}}, u)\right)^2} \\[-0.1cm]
&\hphantom{= i} + \frac{1}{n} \sum_{j=1}^n \left(\left(\hat{F}_1(t \mid A=1, \boldsymbol{Z_j})\right)^2 \exp(2 \hat{\beta}_{kA}) + \left(\hat{F}_1(t \mid A=0, \boldsymbol{Z_j})\right)^2\right) \, \exp(2 \skew{3.5}\hat{\boldsymbol{\beta}}_{\boldsymbol{kZ}}^T \boldsymbol{Z_j}) \int_0^\tau \frac{1}{n}\frac{dN_k(u)}{\left(S^{(0)}(\skew{3.5}\hat{\boldsymbol{\beta}}_{\boldsymbol{k}}, u)\right)^2} \\[-0.1cm]
&\hphantom{= i} + \frac{1}{n} \sum_{j=1}^n \exp(2 \hat{\beta}_{kA} + 2 \skew{3.5}\hat{\boldsymbol{\beta}}_{\boldsymbol{kZ}}^T \boldsymbol{Z_j}) \\[-0.3cm]
&\hphantom{= i + \frac{1}{n} \sum_{j=1}^n i} \cdot \int_0^t \left((1, \boldsymbol{Z}_{\boldsymbol{j}}^T)^T - \boldsymbol{E}(\skew{3.5}\hat{\boldsymbol{\beta}}_{\boldsymbol{k}}, u)\right)^T \hat{\mathbf{\Sigma}}_\mathbf{k}^{-1} \left(\frac{1}{n} \sum_{i=1}^n \tilde{\mathbf{\Sigma}}_\mathbf{ki}\right) \hat{\mathbf{\Sigma}}_\mathbf{k}^{-1} \left((1, \boldsymbol{Z}_{\boldsymbol{j}}^T)^T - \boldsymbol{E}(\skew{3.5}\hat{\boldsymbol{\beta}}_{\boldsymbol{k}}, u)\right) \frac{1}{n} \frac{dN_k(u)}{\left(S^{(0)}(\skew{3.5}\hat{\boldsymbol{\beta}}_{\boldsymbol{k}}, u)\right)^2} \\[-0.1cm]
&\hphantom{= i} + \frac{1}{n} \sum_{j=1}^n \exp(2 \skew{3.5}\hat{\boldsymbol{\beta}}_{\boldsymbol{kZ}}^T \boldsymbol{Z_i}) \\[-0.3cm]
&\hphantom{= i + \frac{1}{n} \sum_{j=1}^n i} \cdot \int_0^t \left((0, \boldsymbol{Z}_{\boldsymbol{i}}^T)^T - \boldsymbol{E}(\skew{3.5}\hat{\boldsymbol{\beta}}_{\boldsymbol{k}}, u)\right)^T \hat{\mathbf{\Sigma}}_\mathbf{k}^{-1} \left(\frac{1}{n} \sum_{i=1}^n \tilde{\mathbf{\Sigma}}_\mathbf{ki}\right) \hat{\mathbf{\Sigma}}_\mathbf{k}^{-1}  \left((0, \boldsymbol{Z}_{\boldsymbol{i}}^T)^T - \boldsymbol{E}(\skew{3.5}\hat{\boldsymbol{\beta}}_{\boldsymbol{k}}, u)\right) \frac{1}{n} \frac{dN_k(u)}{\left(S^{(0)}(\skew{3.5}\hat{\boldsymbol{\beta}}_{\boldsymbol{k}}, u)\right)^2} \\[-0.1cm]
&\hphantom{= i} + \frac{1}{n} \sum_{j=1}^n \left(\hat{F}_1(t \mid A=1, \boldsymbol{Z_j})\right)^2 \exp(2 \hat{\beta}_{kA} + 2 \skew{3.5}\hat{\boldsymbol{\beta}}_{\boldsymbol{kZ}}^T \boldsymbol{Z_j}) \\[-0.3cm]
&\hphantom{= i + \frac{1}{n} \sum_{j=1}^n i} \cdot \int_0^\tau \left((1, \boldsymbol{Z}_{\boldsymbol{j}}^T)^T - \boldsymbol{E}(\skew{3.5}\hat{\boldsymbol{\beta}}_{\boldsymbol{k}}, u)\right)^T \hat{\mathbf{\Sigma}}_\mathbf{k}^{-1} \left(\frac{1}{n} \sum_{i=1}^n \tilde{\mathbf{\Sigma}}_\mathbf{ki}\right) \hat{\mathbf{\Sigma}}_\mathbf{k}^{-1} \left((1, \boldsymbol{Z}_{\boldsymbol{j}}^T)^T - \boldsymbol{E}(\skew{3.5}\hat{\boldsymbol{\beta}}_{\boldsymbol{k}}, u)\right) \frac{1}{n} \frac{dN_k(u)}{\left(S^{(0)}(\skew{3.5}\hat{\boldsymbol{\beta}}_{\boldsymbol{k}}, u)\right)^2} \\[-0.1cm]
&\hphantom{= i} + \frac{1}{n} \sum_{j=1}^n \left(\hat{F}_1(A=0, \boldsymbol{Z_i})\right)^2 \exp(2 \skew{3.5}\hat{\boldsymbol{\beta}}_{\boldsymbol{kZ}}^T \boldsymbol{Z_i}) \\[-0.4cm]
&\hphantom{= i + \frac{1}{n} \sum_{j=1}^n i} \cdot \int_0^\tau \left((0, \boldsymbol{Z}_{\boldsymbol{i}}^T)^T - \boldsymbol{E}(\skew{3.5}\hat{\boldsymbol{\beta}}_{\boldsymbol{k}}, u)\right)^T \hat{\mathbf{\Sigma}}_\mathbf{k}^{-1} \left(\frac{1}{n} \sum_{i=1}^n \tilde{\mathbf{\Sigma}}_\mathbf{ki}\right) \hat{\mathbf{\Sigma}}_\mathbf{k}^{-1}  \left((0, \boldsymbol{Z}_{\boldsymbol{i}}^T)^T - \boldsymbol{E}(\skew{3.5}\hat{\boldsymbol{\beta}}_{\boldsymbol{k}}, u)\right) \frac{1}{n} \frac{dN_k(u)}{\left(S^{(0)}(\skew{3.5}\hat{\boldsymbol{\beta}}_{\boldsymbol{k}}, u)\right)^2}
\end{align*}
with $\smash{\tilde{\mathbf{\Sigma}}_\mathbf{ki} = \int_0^\tau \left((A_i, \boldsymbol{Z}_{\boldsymbol{i}}^T)^T - \boldsymbol{E}(\skew{3.5}\hat{\boldsymbol{\beta}}_{\boldsymbol{k}}, u)\right) \left((A_i, \boldsymbol{Z}_{\boldsymbol{i}}^T)^T - \boldsymbol{E}(\skew{3.5}\hat{\boldsymbol{\beta}}_{\boldsymbol{k}}, u)\right)^T \! dN_{ki}(u)}$, the following inequality holds in probability provided that the conditions in Theorem~\ref{theorem:U_hat} are fulfilled:
\begin{gather*}
\mathbb{E}\left(\left(\sum_{i=1}^n G_i X_{n,i}^{(k)}(t_r) - \sum_{i=1}^n G_i X_{n,i}^{(k)}(t_q)\right)^2 \left(\sum_{i=1}^n G_i X_{n,i}^{(k)}(t_s) - \sum_{i=1}^n G_i X_{n,i}^{(k)}(t_r)\right)^2 \;\biggm|\; \mathscr{F}_\tau\right) \\[-0.1cm]
\leq \left(L_n^{(k)}(t_s) - L_n^{(k)}(t_q)\right)^{3/2} \!\! \cdot O_p(1).
\end{gather*}
\end{lemma}

\begin{proof}[Proof of Theorem~\ref{theorem:U_hat}]
Considering $\smash{\hat{U}_n^{(k)}(\cdot) = \sum_{i=1}^n G_i X_{n,i}^{(k)}(\cdot)}$, the conditions of Lemma~1 from the supplementary material of \citet{dobler2017non} are fulfilled ${\forall k \in \{1, \dots, K\}}$ due to Lemma~\ref{lemma_conditions_X} and the assumptions w.r.t.~the multipliers $\smash{G_i}$. It follows that the finite-dimensional distributions of $\smash{\hat{U}_n^{(k)}}$ converge weakly to zero-mean Gaussian processes with covariance functions $\smash{\xi^{(k)}}$, respectively, in probability (conditional on $\smash{\mathscr{F}_\tau}$). \\
Since $\smash{\tfrac{1}{n} \sum_{i=1}^n \tilde{\mathbf{\Sigma}}_\mathbf{ki}}$ converges to $\smash{\mathbf{\Sigma_k}}$ (cf. proof of Lemma~\ref{lemma_conditions_X}), the function $\smash{L_n^{(k)}}$ from Lemma~\ref{lemma_boundedness_expectation_GX} converges uniformly to 
\begin{align*}
l^{(k)}(t) &= \left(\exp(2\beta_{0kA}) + 2 \exp(\beta_{0kA} + 1)\right) \mathbb{E}_{\boldsymbol{Z}}\left(\exp(2 \boldsymbol{\beta}_{\boldsymbol{0kZ}}^T \boldsymbol{Z})\right) \int_0^t \frac{d \Lambda_{0k}(u)}{s^{(0)}(\boldsymbol{\beta_{0k}}, u)} \\
&\hphantom{= i} + \mathbb{E}_{\boldsymbol{Z}}\left(\left(\left(F_1(t \mid A=1, \boldsymbol{Z})\right)^2 \exp(2\beta_{0kA}) + \left(F_1(t \mid A=0, \boldsymbol{Z})\right)^2\right) \exp(2 \boldsymbol{\beta}_{\boldsymbol{0kZ}}^T \boldsymbol{Z})\right) \int_0^\tau \frac{d \Lambda_{0k}(u)}{s^{(0)}(\boldsymbol{\beta_{0k}}, u)} \\
&\hphantom{= i} + \mathbb{E}_{\boldsymbol{Z}}\left(\exp(2\beta_{0kA} + 2 \boldsymbol{\beta}_{\boldsymbol{0kZ}}^T \boldsymbol{Z}) \int_0^t \left((1, \boldsymbol{Z}^T)^T - \boldsymbol{e}(\boldsymbol{\beta_{0k}}, u)\right)^T \mathbf{\Sigma_k}^{\!\!\!\!-1} \left((1, \boldsymbol{Z}^T)^T - \boldsymbol{e}(\boldsymbol{\beta_{0k}}, u)\right) \frac{d \Lambda_{0k}(u)}{s^{(0)}(\boldsymbol{\beta_{0k}}, u)} \right) \\
&\hphantom{= m} + \mathbb{E}_{\boldsymbol{Z}}\left(\exp(2 \boldsymbol{\beta}_{\boldsymbol{0kZ}}^T \boldsymbol{Z}) \int_0^t \left((0, \boldsymbol{Z}^T)^T - \boldsymbol{e}(\boldsymbol{\beta_{0k}}, u)\right)^T \mathbf{\Sigma_k}^{\!\!\!\!-1} \left((0, \boldsymbol{Z}^T)^T - \boldsymbol{e}(\boldsymbol{\beta_{0k}}, u)\right) \frac{d \Lambda_{0k}(u)}{s^{(0)}(\boldsymbol{\beta_{0k}}, u)} \right) \\
&\hphantom{= i} + \mathbb{E}_{\boldsymbol{Z}}\left(\left(F_1(t \mid A=1, \boldsymbol{Z})\right)^2 \exp(2\beta_{0kA} + 2 \boldsymbol{\beta}_{\boldsymbol{0kZ}}^T \boldsymbol{Z}) \int_0^\tau \left((1, \boldsymbol{Z}^T)^T - \boldsymbol{e}(\boldsymbol{\beta_{0k}}, u)\right)^T \mathbf{\Sigma_k}^{\!\!\!\!-1} \left((1, \boldsymbol{Z}^T)^T - \boldsymbol{e}(\boldsymbol{\beta_{0k}}, u)\right) \right. \\[-0.2cm]
&\hphantom{= i + \mathbb{E}_{\boldsymbol{Z}}\left(\left(F_1(t \mid A=1, \boldsymbol{Z})\right)^2 \exp(2\beta_{0kA} + 2 \boldsymbol{\beta}_{\boldsymbol{0kZ}}^T \boldsymbol{Z}) \int_0^\tau i\right.} \left. \cdot \frac{d \Lambda_{0k}(u)}{s^{(0)}(\boldsymbol{\beta_{0k}}, u)} \right) \\
&\hphantom{= i} + \mathbb{E}_{\boldsymbol{Z}}\left(\left(F_1(t \mid A=0, \boldsymbol{Z})\right)^2 \exp(2 \boldsymbol{\beta}_{\boldsymbol{0kZ}}^T \boldsymbol{Z}) \int_0^\tau \left((0, \boldsymbol{Z}^T)^T - \boldsymbol{e}(\boldsymbol{\beta_{0k}}, u)\right)^T \mathbf{\Sigma_k}^{\!\!\!\!-1} \left((0, \boldsymbol{Z}^T)^T - \boldsymbol{e}(\boldsymbol{\beta_{0k}}, u)\right) \right. \\[-0.2cm]
&\hphantom{= i + \mathbb{E}_{\boldsymbol{Z}}\left(\left(F_1(t \mid A=0, \boldsymbol{Z})\right)^2 \exp(2 \boldsymbol{\beta}_{\boldsymbol{0kZ}}^T \boldsymbol{Z}) \int_0^\tau i\right.} \left. \cdot \frac{d \Lambda_{0k}(u)}{s^{(0)}(\boldsymbol{\beta_{0k}}, u)} \right)
\end{align*}
on ${[0, \tau)}$ as a consequence of the martingale central limit theorem. The conditional tightness of $\smash{\hat{U}_n^{(k)}}$ can now be shown along the lines of the proof of Theorem~3.1 in \citet{dobler2014bootstrapping}. We apply the subsequence principle for convergence in probability \citep[cf.][]{beyersmann2013weak}: For every subsequence, there is another subsequence such that for almost every (fixed) $\smash{\omega \in \Omega_1\!\times\!\Omega_2}$, we find ${n_0 \in \mathbb{N}}$, a constant ${\gamma > 0}$ and a sequence of non-decreasing, continuous functions $\smash{l_n^{(k)}}$ that converges uniformly to $\smash{l^{(k)}}$, such that
\[\mathbb{E}\left(\left(\sum_{i=1}^n G_i X_{n,i}^{(k)}(t_r) - \sum_{i=1}^n G_i X_{n,i}^{(k)}(t_q)\right)^2 \left(\sum_{i=1}^n G_i X_{n,i}^{(k)}(t_s) - \sum_{i=1}^n G_i X_{n,i}^{(k)}(t_r)\right)^2 \mid \mathscr{F}_\tau\right) \leq \gamma \left(l_n^{(k)}(t_s) - l_n^{(k)}(t_q)\right)^{3/2}\]
if $\smash{n \geq n_0}$. (Here, $\smash{n_0}$ and $\gamma$ do not depend on $\smash{0 \leq t_q \leq t_r \leq t_s \leq \tau}$.) The conditional tightness follows by extending Theorem~13.5 in \citet{billingsley1999convergence} pointwise along subsequences \citep[cf.][]{dobler2014bootstrapping}. Eventually, this proves the conditional convergence in distribution of $\smash{\hat{U}_n^{(k)}}$ in probability for each ${k \in \{1, \dots, K\}}$. \\
The assertion of Theorem~\ref{theorem:U_hat} follows by noting that the processes $\smash{U_n^{(k)}}$ and $\smash{\hat{U}_n^{(k')}}$ are independent for ${k \neq k'}$ given the data because we consider competing events, i.e.~$\smash{dN_{ki}}$ and $\smash{dN_{k'i}}$ cannot jump both.
\end{proof}

\section{Discussion \label{sec:discussion}}

Estimating a causal effect in time-to-event data subject to competing risks presents a number of challenges. In this manuscript, we focused on the average treatment effect defined as the difference between the $t$-year absolute risks for the event of interest. While we based our estimation on cause-specific Cox models, our estimand of interest was the causal risk difference, i.e.~a contrast of the cumulative incidence functions. Thus, we do not interpret the hazard ratio in a causal way, see \citet{hernan2010the,aalen2015does,martinussen2013on} for detailed discussions on the drawbacks of the hazard ratio in a causal context. As \citet{martinussen2023} point out, this estimand only captures the total effect of the treatment on the event of interest, while a distinction in terms of direct and indirect effects is not possible, see also \citet{young2020a}. \citet{martinussen2023} therefore propose an estimator based on the efficient influence function and use the nonparametric bootstrap to estimate its variance. Extensions of the wild bootstrap to this situation merit further research, as the classical bootstrap has been shown to perform insufficiently in certain situations \citep{singh1981on, friedrich2017permuting, niessl2023statistical, ruehl2022general}.

A related aspect concerns modeling of the association between the covariates and the outcome. In our work, we focused on Cox proportional hazards models, but did not go into aspects such as variable selection or model misspecification. The latter is covered to some extent by \citet{ozenne2020on} considering the classical bootstrap. Recently, \citet{vansteelandt2022} proposed an approach that allows for more flexible modeling of the association between covariates and an outcome. Integration of this so-called assumption-lean Cox regression into our resampling framework is part of future research. Alternatively, other regression models such as Aalen's additive hazards model or a Cox-Aalen model might be used, however, the proofs presented here need to be adapted accordingly. Exploiting the properties of the martingale residuals underlying these models will be helpful and, in addition, facilitate the integration of e.g. left-truncation.

\begin{acknowledgments}
Support by the DFG (Grant FR 4121/2-1) is gratefully acknowledged.
The authors further thank Dennis Dobler for his constructive ideas.
\end{acknowledgments}

\appendix \small

\begin{proof}[Proof of Lemma \ref{lemma_conditions_X}]
Because $\smash{N_{ki}}$ jumps at most once, 
\[\max_{1 \leq i \leq n} \left| \int_0^{t_r} \hat{\tilde{H}}_{k1}(u, t_r) \frac{dN_{ki}(u)}{\sqrt{n} \!\cdot\! S^{(0)}(\skew{3.5}\hat{\boldsymbol{\beta}}_{\boldsymbol{k}}, u)}\right| < \left(\exp(\hat{\beta}_{kA}) + 1\right) \max_{1 \leq i \leq n} \exp(\skew{3.5}\hat{\boldsymbol{\beta}}_{\boldsymbol{kZ}}^T \boldsymbol{Z_i}) \frac{1}{\sqrt{n} \!\cdot\! \inf_{u \in [0, t_r]} S^{(0)}(\skew{3.5}\hat{\boldsymbol{\beta}}_{\boldsymbol{k}}, u)}.\]
Recall that on ${\mathcal{B} \times [0, \tau]}$, $\smash{S^{(0)}}$ converges uniformly to $\smash{s^{(0)}}$, which is bounded away from zero, and that $\smash{\skew{3.5}\hat{\boldsymbol{\beta}}_{\boldsymbol{k}}}$ is strongly consistent. For that reason, the expression above converges to 0 $\smash{\forall t_r \in \{t_1, \dots, t_l\}}$ almost surely as ${n \to \infty}$. \\
In addition,
\begin{align*}
\max_{1 \leq i \leq n} &\left| \int_0^\tau \frac{1}{\sqrt{n}} \, \smash{\left(\boldsymbol{\hat{\tilde{H}}_{k2}}(t_r)\right)^T} \! \hat{\mathbf{\Sigma}}_\mathbf{k}^{-1} \left((A_i, \boldsymbol{Z}_{\boldsymbol{i}}^T)^T - \boldsymbol{E}(\skew{3.5}\hat{\boldsymbol{\beta}}_{\boldsymbol{k}}, u)\right) \, dN_{ki}(u) \right| \\
\leq & \; \frac{1}{\sqrt{n}} \max_{1 \leq i \leq n} \left\{ \frac{1}{n} \sum_{j_1=1}^n \frac{1}{n} \sum_{j_2=1}^n \int_0^{t_r} \left(\exp(\hat{\beta}_{kA}) \!\cdot\! \left|\left((1, \boldsymbol{Z}_{\boldsymbol{j_1}}^T) - \left(\boldsymbol{E}(\skew{3.5}\hat{\boldsymbol{\beta}}_{\boldsymbol{k}}, v)\right)^T\right) \hat{\mathbf{\Sigma}}_\mathbf{k}^{-1} \! \int_0^\tau \left((A_i, \boldsymbol{Z}_{\boldsymbol{i}}^T)^T - \boldsymbol{E}(\skew{3.5}\hat{\boldsymbol{\beta}}_{\boldsymbol{k}}, u)\right) \, dN_{ki}(u) \right| \right. \right. \\[-0.2cm]
& \; \hphantom{\frac{1}{\sqrt{n}} \max_{1 \leq i \leq n} \left\{ \frac{1}{n} \sum_{j_1=1}^n \frac{1}{n} \sum_{j_2=1}^n \int_0^{t_r} \left( i \right. \right.} \left.\left. + \left| \left((0, \boldsymbol{Z}_{\boldsymbol{j_1}}^T) - \left(\boldsymbol{E}(\skew{3.5}\hat{\boldsymbol{\beta}}_{\boldsymbol{k}}, v)\right)^T\right) \hat{\mathbf{\Sigma}}_\mathbf{k}^{-1} \! \int_0^\tau \left((A_i, \boldsymbol{Z}_{\boldsymbol{i}}^T)^T - \boldsymbol{E}(\skew{3.5}\hat{\boldsymbol{\beta}}_{\boldsymbol{k}}, u)\right) \, dN_{ki}(u) \right| \right) \right. \\[-0.1cm]
& \; \hphantom{\frac{1}{\sqrt{n}} \max_{1 \leq i \leq n} \left\{ \frac{1}{n} \sum_{j_1=1}^n \frac{1}{n} \sum_{j_2=1}^n \int_0^{t_r} i \right.} \left. \cdot \frac{\exp(\skew{3.5}\hat{\boldsymbol{\beta}}_{\boldsymbol{kZ}}^T \boldsymbol{Z_{j_1}})}{S^{(0)}(\skew{3.5}\hat{\boldsymbol{\beta}}_{\boldsymbol{k}}, v)} \, dN_{k j_2}(v) \right\} \\
<& \; \frac{\max\{\exp(\hat{\beta}_{kA}), 1\} \, \max_{1 \leq j_1 \leq n} \exp(\skew{3.5}\hat{\boldsymbol{\beta}}_{\boldsymbol{kZ}}^T \boldsymbol{Z_{j_1}})}{\sqrt{n} \!\cdot\! \inf_{v \in [0, t_r]} S^{(0)}(\skew{3.5}\hat{\boldsymbol{\beta}}_{\boldsymbol{k}}, v)} \\[-0.1cm]
& \; \hphantom{i} \cdot \left(\max_{1 \leq i, j_1 \leq n} \, \sup_{u \in [0, \tau], \, v \in [0, t_r]} \left| \left((1, \boldsymbol{Z}_{\boldsymbol{j_1}}^T) - \left(\boldsymbol{E}(\skew{3.5}\hat{\boldsymbol{\beta}}_{\boldsymbol{k}}, v)\right)^T\right) \hat{\mathbf{\Sigma}}_\mathbf{k}^{-1} \! \left((A_i, \boldsymbol{Z}_{\boldsymbol{i}}^T)^T - \boldsymbol{E}(\skew{3.5}\hat{\boldsymbol{\beta}}_{\boldsymbol{k}}, u)\right)\right| \right. \\[-0.1cm]
& \, \hphantom{i \cdot \left(i\right.} \left. + \max_{1 \leq i, j_1 \leq n} \, \sup_{u \in [0, \tau], \, v \in [0, t_r]} \left| \left((0, \boldsymbol{Z}_{\boldsymbol{j_1}}^T) - \left(\boldsymbol{E}(\skew{3.5}\hat{\boldsymbol{\beta}}_{\boldsymbol{k}}, v)\right)^T\right) \hat{\mathbf{\Sigma}}_\mathbf{k}^{-1} \! \left((A_i, \boldsymbol{Z}_{\boldsymbol{i}}^T)^T - \boldsymbol{E}(\skew{3.5}\hat{\boldsymbol{\beta}}_{\boldsymbol{k}}, u)\right)\right|\right).
\end{align*}
Using the previous considerations and the fact that $\smash{s^{(1)}}$ and $\smash{s^{(2)}}$ are bounded on ${\mathcal{B} \times [0, \tau]}$ \citep[section~8.4]{fleming2005counting}, we can also conclude that the above maximum vanishes $\smash{\forall t_r \in \{t_1, \dots, t_l\}}$ as $n$ tends to $\infty$, which implies condition (i).

Moreover, for time points $\smash{t_r}$ and $\smash{t_s}$ with $\smash{0 \leq t_r \leq t_s \leq \tau}$,
\begin{equation} \label{eq:xi}
\begin{aligned} 
\sum_{i=1}^n X_{n,i}^{(k)}(t_r) X_{n,i}^{(k)}(t_s) &= \frac{1}{n} \sum_{i=1}^n \int_0^{t_r} \hat{\tilde{H}}_{k1}(u, t_r) \hat{\tilde{H}}_{k1}(u, t_s) \frac{dN_{ki}(u)}{\left(S^{(0)}(\skew{3.5}\hat{\boldsymbol{\beta}}_{\boldsymbol{k}}, u)\right)^2} \\[-0.2cm]
&\hphantom{= i} + \frac{1}{n} \sum_{i=1}^n \int_0^{t_r} \hat{\tilde{H}}_{k1}(u, t_r) \left(\boldsymbol{\hat{\tilde{H}}_{k2}}(t_s)\right)^T \! \hat{\mathbf{\Sigma}}_\mathbf{k}^{-1} \left((A_i, \boldsymbol{Z}_{\boldsymbol{i}}^T)^T - \boldsymbol{E}(\skew{3.5}\hat{\boldsymbol{\beta}}_{\boldsymbol{k}}, u)\right) \frac{dN_{ki}(u)}{S^{(0)}(\skew{3.5}\hat{\boldsymbol{\beta}}_{\boldsymbol{k}}, u)} \\[-0.1cm]
&\hphantom{= i} + \frac{1}{n} \sum_{i=1}^n \int_0^{t_s} \left(\boldsymbol{\hat{\tilde{H}}_{k2}}(t_r)\right)^T \! \hat{\mathbf{\Sigma}}_\mathbf{k}^{-1} \left((A_i, \boldsymbol{Z}_{\boldsymbol{i}}^T)^T - \boldsymbol{E}(\skew{3.5}\hat{\boldsymbol{\beta}}_{\boldsymbol{k}}, u)\right) \hat{\tilde{H}}_{k1}(u, t_s) \frac{dN_{ki}(u)}{S^{(0)}(\skew{3.5}\hat{\boldsymbol{\beta}}_{\boldsymbol{k}}, u)} \\[-0.1cm]
&\hphantom{= i} + \frac{1}{n} \sum_{i=1}^n \int_0^\tau \left(\boldsymbol{\hat{\tilde{H}}_{k2}}(t_r)\right)^T \! \hat{\mathbf{\Sigma}}_\mathbf{k}^{-1} \left((A_i, \boldsymbol{Z}_{\boldsymbol{i}}^T)^T - \boldsymbol{E}(\skew{3.5}\hat{\boldsymbol{\beta}}_{\boldsymbol{k}}, u)\right) \\[-0.3cm]
&\hphantom{= i + \frac{1}{n} \sum_{i=1}^n \int_0^\tau i} \cdot \left(\boldsymbol{\hat{\tilde{H}}_{k2}}(t_s)\right)^T \! \hat{\mathbf{\Sigma}}_\mathbf{k}^{-1} \left((A_i, \boldsymbol{Z}_{\boldsymbol{i}}^T)^T - \boldsymbol{E}(\skew{3.5}\hat{\boldsymbol{\beta}}_{\boldsymbol{k}}, u)\right) \, dN_{ki}(u),
\end{aligned}
\end{equation}
as $\smash{N_{ki}}$ is a one-jump process. The first term of Equation (\ref{eq:xi}) equals 
\[\frac{1}{n} \sum_{i=1}^n \int_0^{t_r} \hat{\tilde{H}}_{k1}(u, t_r) \hat{\tilde{H}}_{k1}(u, t_s) \frac{dM_{ki}(u)}{\left(S^{(0)}(\skew{3.5}\hat{\boldsymbol{\beta}}_{\boldsymbol{k}}, u)\right)^2} + \int_0^{t_r} \hat{\tilde{H}}_{k1}(u, t_r) \hat{\tilde{H}}_{k1}(u, t_s) \frac{d \Lambda_{0k}(u)}{S^{(0)}(\skew{3.5}\hat{\boldsymbol{\beta}}_{\boldsymbol{k}}, u)}.\]
Due to the strong consistency of $\smash{\skew{3.5}\hat{\boldsymbol{\beta}}_{\boldsymbol{k}}}$ and $\smash{\hat{\Lambda}_{0k}}$, $\smash{\hat{\tilde{H}}_{k1}}$ is uniformly consistent, and so is $\smash{S^{(0)}}$ on ${\mathcal{B} \times [0, \tau]}$ (with estimand $\smash{s^{(0)}}$, which is bounded away from zero). It follows by application of the martingale central limit theorem that the first summand of the expression above converges to zero as ${n \to \infty}$. Using the same arguments on the remaining terms in Equation (\ref{eq:xi}), we obtain
\begin{align*}
\sum_{i=1}^n X_{n,i}^{(k)}(t_r) X_{n,i}^{(k)}(t_s) \overset{P}{\longrightarrow} &\int_0^{t_r} \tilde{H}_{k1}(u, t_r) \tilde{H}_{k1}(u, t_s) \frac{d \Lambda_{0k}(u)}{s^{(0)}(\boldsymbol{\beta_{0k}},u)} \\[-0.2cm]
&\hphantom{i} + \int_0^{t_r} \tilde{H}_{k1}(u, t_r) \left(\boldsymbol{\tilde{H}_{k2}}(t_s)\right)^T \! \mathbf{\Sigma_k}^{\!\!\!\!-1} \left(\boldsymbol{s^{(1)}}(\boldsymbol{\beta_{0k}}, u) - \boldsymbol{e}(\boldsymbol{\beta_{0k}}, u) s^{(0)}(\boldsymbol{\beta_{0k}}, u)\right) \frac{d \Lambda_{0k}(u)}{s^{(0)}(\boldsymbol{\beta_{0k}}, u)} \\[-0.1cm]
&\hphantom{i} + \int_0^{t_s} \left(\boldsymbol{\tilde{H}_{k2}}(t_r)\right)^T \! \mathbf{\Sigma_k}^{\!\!\!\!-1} \left(\boldsymbol{s^{(1)}}(\boldsymbol{\beta_{0k}}, u) - \boldsymbol{e}(\boldsymbol{\beta_{0k}}, u) s^{(0)}(\boldsymbol{\beta_{0k}}, u)\right) \tilde{H}_{k1}(u, t_s) \frac{d \Lambda_{0k}(u)}{s^{(0)}(\boldsymbol{\beta_{0k}}, u)} \\[-0.1cm]
&\hphantom{i} + \left(\boldsymbol{\tilde{H}_{k2}}(t_r)\right)^T \! \mathbf{\Sigma_k}^{\!\!\!\!-1} \left(\int_0^\tau \left(\mathbf{s^{(2)}}(\boldsymbol{\beta_{0k}}, u) - \boldsymbol{s^{(1)}}(\boldsymbol{\beta_{0k}}, u) \left(\boldsymbol{e}(\boldsymbol{\beta_{0k}}, u)\right)^T\right) d \Lambda_{0k}(u)\right) \left(\mathbf{\Sigma_k}^{\!\!\!\!-1}\right)^T \! \boldsymbol{\tilde{H}_{k2}}(t_s) \\
&= \int_0^{t_r} \tilde{H}_{k1}(u, t_r) \tilde{H}_{k1}(u, t_s) \frac{d \Lambda_{0k}(u)}{s^{(0)}(\boldsymbol{\beta_{0k}},u)} + \left(\boldsymbol{\tilde{H}_{k2}}(t_r)\right)^T \! \mathbf{\Sigma_k}^{\!\!\!\!-1} \mathbf{\Sigma_k} \mathbf{\Sigma_k}^{\!\!\!\!-1} \, \boldsymbol{\tilde{H}_{k2}}(t_s),
\end{align*}
and thus, condition (ii) follows.
\end{proof}

\begin{proof}[Proof of Lemma \ref{lemma_boundedness_X}]
\begin{align*}
\sqrt{n} \max_{1 \leq i \leq n} \left|X_{n,i}^{(k)}(t_s) - X_{n,i}^{(k)}(t_r) \right| &\leq \max_{1 \leq i \leq n} \left\{\int_0^{t_s} \left|\hat{\tilde{H}}_{k1}(u, t_s) - \mathbbm{1}\{u \leq t_r\} \!\cdot\! \hat{\tilde{H}}_{k1}(u, t_r) \right| \frac{dN_{ki}(u)}{S^{(0)}(\skew{3.5}\hat{\boldsymbol{\beta}}_{\boldsymbol{k}}, u)} \right. \\[-0.1cm]
&\hphantom{\leq \max_{1 \leq i \leq n} \left\{ i \right.} \left. + \int_0^\tau \left|\smash{\left(\boldsymbol{\hat{\tilde{H}}_{k2}}(t_s) - \boldsymbol{\hat{\tilde{H}}_{k2}}(t_r)\right)^T} \! \hat{\mathbf{\Sigma}}_\mathbf{k}^{-1} \left((A_i, \boldsymbol{Z}_{\boldsymbol{i}}^T)^T - \boldsymbol{E}(\skew{3.5}\hat{\boldsymbol{\beta}}_{\boldsymbol{k}}, u)\right)\right| dN_{ki}(u) \right\} \\
&< \frac{2 \left(\exp(\smash{\hat{\beta}_{kA}}) + 1\right) \max_{1 \leq i \leq n} \exp(\skew{3.5}\hat{\boldsymbol{\beta}}_{\boldsymbol{kZ}}^T \boldsymbol{Z_i})}{\inf_{u \in [0, \tau]} S^{(0)}(\skew{3.5}\hat{\boldsymbol{\beta}}_{\boldsymbol{k}}, u)} \\
&\hphantom{< i} + \max_{1 \leq i \leq n} \, \sup_{u, t_r, t_s \in [0, \tau]} \left|\smash{\left(\boldsymbol{\hat{\tilde{H}}_{k2}}(t_s) - \boldsymbol{\hat{\tilde{H}}_{k2}}(t_r)\right)^T} \! \hat{\mathbf{\Sigma}}_\mathbf{k}^{-1} \left((A_i, \boldsymbol{Z}_{\boldsymbol{i}}^T)^T - \boldsymbol{E}(\skew{3.5}\hat{\boldsymbol{\beta}}_{\boldsymbol{k}}, u)\right)\right|,
\end{align*}
i.e.~$\smash{\sqrt{n} \max_{1 \leq i \leq n} \left|X_{n,i}^{(k)}(t_s) - X_{n,i}^{(k)}(t_r) \right| \in O_p(1)}$ (cf.~the proof of Lemma~\ref{lemma_conditions_X}).
\end{proof}

\begin{proof}[Proof of Lemma \ref{lemma_boundedness_expectation_GX}]
Using condition (iv) of Theorem~\ref{theorem:U_hat}, one can show that the expectation in Lemma~\ref{lemma_boundedness_expectation_GX} has the upper bound
\begin{gather*}
\max_{1 \leq i \leq n} \mathbb{E}(G_i^4 \mid \mathscr{F}_\tau) \sum_{i=1}^n \left(X_{n,i}^{(k)}(t_r) - X_{n,i}^{(k)}(t_q)\right)^2 \left(X_{n,i}^{(k)}(t_s) - X_{n,i}^{(k)}(t_r)\right)^2 \\
\begin{multlined}
+ 2 \max_{1 \leq i_1 \leq n} \left|\mathbb{E}\left(G_{i_1}^3 \mid \mathscr{F}_\tau\right)\right| \max_{1 \leq i_2 \leq n} \left|\mathbb{E}\left(G_{i_2} \mid \mathscr{F}_\tau\right)\right| \sum_{i_1=1}^n \left(X_{n, i_1}^{(k)} (t_r) - X_{n, i_1}^{(k)} (t_q)\right)^2 \left|X_{n, i_1}^{(k)}(t_s) - X_{n, i_1}^{(k)}(t_r)\right| \\[-0.3cm]
\cdot \sum_{i_2=1}^n \left|X_{n, i_2}^{(k)}(t_s) - X_{n, i_2}^{(k)}(t_r)\right|
\end{multlined} \\
\begin{multlined}
+ 2 \max_{1 \leq i_1 \leq n} \left|\mathbb{E}\left(G_{i_1} \mid \mathscr{F}_\tau\right)\right| \max_{1 \leq i_2 \leq n} \left|\mathbb{E}\left(G_{i_2}^3 \mid \mathscr{F}_\tau\right)\right| \sum_{i_1=1}^n \left|X_{n, i_1}^{(k)} (t_r) - X_{n, i_1}^{(k)} (t_q)\right| \\[-0.3cm]
\cdot \sum_{i_2=1}^n \left|X_{n, i_1}^{(k)}(t_r) - X_{n, i_1}^{(k)}(t_q)\right|  \left(X_{n, i_2}^{(k)}(t_s) - X_{n, i_2}^{(k)}(t_r)\right)^2
\end{multlined} \\
+ \max_{1 \leq i \leq n} \left(\mathbb{E}\left(G_i^2 \mid \mathscr{F}_\tau\right)\right)^2 \sum_{i_1=1}^n \left(X_{n, i_1}^{(k)}(t_r) - X_{n, i_1}^{(k)}(t_q)\right)^2 \sum_{i_2=1}^n \left(X_{n, i_2}^{(k)}(t_s) - X_{n, i_2}^{(k)}(t_r)\right)^2 \\
+ 2 \max_{1 \leq i \leq n} \left(\mathbb{E}\left(G_i^2 \mid \mathscr{F}_\tau\right)\right)^2 \left(\sum_{i=1}^n \left|X_{n, i}^{(k)}(t_r) - X_{n, i}^{(k)}(t_q)\right| \left|X_{n, i}^{(k)}(t_s) - X_{n, i}^{(k)}(t_r)\right|\right)^2 \\
+ \max_{1 \leq i_1 \leq n} \mathbb{E}\left(G_{i_1}^2 \mid \mathscr{F}_\tau\right) \max_{1 \leq i \leq n} \left(\mathbb{E}\left(G_i \mid \mathscr{F}_\tau\right)\right)^2 \sum_{i_1=1}^n \left(X_{n, i_1}^{(k)}(t_r) - X_{n, i_1}^{(k)}(t_q)\right)^2 \left(\sum_{i=1}^n \left|X_{n, i}^{(k)}(t_s) - X_{n, i}^{(k)}(t_r)\right|\right)^2 \\
\begin{multlined}
+ 4 \max_{1 \leq i_2 \leq n} \mathbb{E}\left(G_{i_2}^2 \mid \mathscr{F}_\tau\right) \max_{1 \leq i \leq n} \left(\mathbb{E}\left(G_i^2 \mid \mathscr{F}_\tau\right)\right)^2 \sum_{i_1=1}^n \left|X_{n, i_1}^{(k)}(t_r) - X_{n, i_1}^{(k)}(t_q)\right| \\[-0.3cm]
\cdot \sum_{i_2=1}^n \left|X_{n, i_2}^{(k)}(t_r) - X_{n, i_2}^{(k)}(t_q)\right| \left|X_{n, i_2}^{(k)}(t_s) - X_{n, i_2}^{(k)}(t_r)\right| \sum_{i_3=1}^n \left|X_{n, i_3}^{(k)}(t_s) - X_{n, i_3}^{(k)}(t_r)\right|
\end{multlined} \\
+ \max_{1 \leq i \leq n} \left(\mathbb{E}\left(G_i \mid \mathscr{F}_\tau\right)\right)^2 \max_{1 \leq i_3 \leq n} \mathbb{E}\left(G_{i_3}^2 \mid \mathscr{F}_\tau\right) \left(\sum_{i=1}^n \left|X_{n, i}^{(k)}(t_r) - X_{n, i}^{(k)}(t_q)\right|\right)^2 \sum_{i_3=1}^n \left(X_{n, i_3}^{(k)}(t_s) - X_{n, i_3}^{(k)}(t_r)\right)^2 \\
+ \max_{1 \leq i \leq n} \left(\mathbb{E}\left(G_i \mid \mathscr{F}_\tau\right)\right)^4 \left(\sum_{i_1=1}^n \left|X_{n, i_1}^{(k)}(t_r) - X_{n, i_1}^{(k)}(t_q)\right|\right)^2 \left(\sum_{i_2=1}^n \left|X_{n, i_2}^{(k)}(t_s) - X_{n, i_2}^{(k)}(t_r)\right|\right)^2,
\end{gather*}
which we denote by $\refstepcounter{equation} (\theequation) \label{eq:lemma4_1}$. According to the proof of Lemma~\ref{lemma_conditions_X}, the first term can (informally) be expressed as 
\[\max_{1 \leq i \leq n} \mathbb{E}(G_i^4 \mid \mathscr{F}_\tau) \, \frac{1}{n^2} \sum_{i=1}^n \left(\int_0^{t_r} dN_{ki}(u) \cdot O_p(1) + \int_0^\tau dN_{ki}(u) \cdot O_p(1)\right)^2 \left(\int_0^{t_s} dN_{ki}(u) \cdot O_p(1) + \int_0^\tau dN_ki(u) \cdot O_p(1)\right)^2,\]
which may be further reduced to $\smash{\max_{1 \leq i \leq n} \mathbb{E}(G_i^4 \mid \mathscr{F}_\tau) \, \frac{1}{n} \cdot O_p(1)}$, as $\smash{N_{ki}}$ is a one-jump process. The term at hand is therefore negligible due to condition (iii) of the theorem. \\
Furthermore, the second and third summands in (\ref{eq:lemma4_1}) have the upper bound
\[\max_{1 \leq i_1 \leq n} \left|\mathbb{E}\left(G_{i_1}^3 \mid \mathscr{F}_\tau\right)\right| \max_{1 \leq i_2 \leq n} \left|\mathbb{E}\left(G_{i_2} \mid \mathscr{F}_\tau\right)\right| \max_{(t_o, t_p) \in \{(t_q,t_r), (t_r,t_s)\}} \left(\sum_{i=1}^n \left(X_{n,i}^{(k)}(t_p) - X_{n,i}^{(k)}(t_o)\right)^2\right)^{3/2} \! \sqrt{n} \; O_p(n^{-1/2})\]
as a consequence of the Cauchy-Schwarz inequality and Lemma~\ref{lemma_boundedness_X}. With condition (i) as well as a combination of the Jensen inequality and condition (iii), we eventually obtain the representation
\[\max_{(t_o, t_p) \in \{(t_q,t_r), (t_r,t_s)\}} \left(\sum_{i=1}^n \left(X_{n,i}^{(k)}(t_p) - X_{n,i}^{(k)}(t_o)\right)^2\right)^{3/2} O_p(1).\]
This expression turns out to be a general upper bound for the expectation in Lemma~\ref{lemma_boundedness_expectation_GX} by application of similar considerations, involving the Cauchy-Schwarz inequality, Lemma~\ref{lemma_boundedness_X} and the conditions of Theorem~\ref{theorem:U_hat}, to the remaining terms in (\ref{eq:lemma4_1}). Note that the $\smash{O_p(1)}$ term does not depend on $\smash{t_q, t_r, t_s}$!

For $\smash{(t_o, t_p) \in \{(t_q, t_r), (t_r, t_s)\}}$, it thus remains to show that 
\[\sum_{i=1}^n \left(X_{n,i}^{(k)}(t_p) - X_{n,i}^{(k)}(t_o)\right)^2 \leq \left(L_n^{(k)}(t_s) - L_n^{(k)}(t_q)\right) O_p(1).\]
The inequality $\smash{(a+b)^2 \leq 2a^2 + 2b^2}$, $a, b \in \mathbb{R}$ suggests that
\begin{equation} \label{eq:lemma4_2}
\begin{aligned}
\sum_{i=1}^n \left(X_{n,i}^{(k)}(t_p) - X_{n,i}^{(k)}(t_o)\right)^2 &\leq \frac{2}{n} \sum_{i=1}^n \left(\left(\int_0^{t_p} \hat{\tilde{H}}_{k1}(u, t_p) \frac{dN_{ki}(u)}{S^{(0)}(\skew{3.5}\hat{\boldsymbol{\beta}}_{\boldsymbol{k}}, u)} - \int_0^{t_o} \hat{\tilde{H}}_{k1}(u, t_0) \frac{dN_{ki}(u)}{S^{(0)}(\skew{3.5}\hat{\boldsymbol{\beta}}_{\boldsymbol{k}}, u)}\right)^2 \right. \\[-0.1cm]
&\hphantom{\leq \frac{2}{n} \sum_{i=1}^n \left( i \right.} + \left.\left(\int_0^\tau \left(\boldsymbol{\hat{\tilde{H}}_{k2}}(t_p) - \boldsymbol{\hat{\tilde{H}}_{k2}}(t_o)\right)^T \! \hat{\mathbf{\Sigma}}_\mathbf{k}^{-1} \left((A_i, \boldsymbol{Z}_{\boldsymbol{i}}^T)^T - \boldsymbol{E}(\skew{3.5}\hat{\boldsymbol{\beta}}_{\boldsymbol{k}}, u)\right) dN_{ki}(u)\right)^2\right),
\end{aligned}
\end{equation}
Due to the definition of $\smash{\hat{\tilde{H}}_{k1}}$, the first summand in (\ref{eq:lemma4_2}) has upper bound
\begin{equation} \label{eq:lemma4_3}
\begin{aligned}
&\frac{2}{n} \sum_{i=1}^n \left(\frac{1}{n} \sum_{j=1}^n  \left(2 \exp(\hat{\beta}_{kA}) + 2\right) \exp(\skew{3.5}\hat{\boldsymbol{\beta}}_{\boldsymbol{kZ}}^T \boldsymbol{Z_j}) \int_{t_o}^{t_p} \frac{dN_{ki}(u)}{S^{(0)}(\skew{3.5}\hat{\boldsymbol{\beta}}_{\boldsymbol{k}}, u)} +  \frac{1}{n} \sum_{j=1}^n  \left(\exp(\hat{\beta}_{kA}) + 1\right) \exp(\skew{3.5}\hat{\boldsymbol{\beta}}_{\boldsymbol{kZ}}^T \boldsymbol{Z_j}) \int_{t_o}^{t_p} \frac{dN_{ki}(u)}{S^{(0)}(\skew{3.5}\hat{\boldsymbol{\beta}}_{\boldsymbol{k}}, u)} \right. \\[-0.1cm]
&\hphantom{\frac{2}{n} \sum_{i=1}^n \left( i \right.} \left. + \frac{1}{n} \sum_{j=1}^n \left(\left(\hat{F}_1(t_p \mid A=1, \boldsymbol{Z_j}) - \hat{F}_1(t_o \mid A=1, \boldsymbol{Z_j})\right) \exp(\hat{\beta}_{kA}) + \left(\hat{F}_1(t_p \mid A=0, \boldsymbol{Z_j}) - \hat{F}_1(t_o \mid A=0, \boldsymbol{Z_j})\right) \right) \right. \\[-0.3cm]
&\hphantom{\frac{2}{n} \sum_{i=1}^n \left( i + \frac{1}{n} \sum_{j=1}^n i \right.} \left. \cdot \exp(\skew{3.5}\hat{\boldsymbol{\beta}}_{\boldsymbol{kZ}}^T \boldsymbol{Z_j}) \int_0^{t_o} \frac{dN_{ki}(u)}{S^{(0)}(\skew{3.5}\hat{\boldsymbol{\beta}}_{\boldsymbol{k}}, u)}\right)^2 \\
&\leq \frac{2}{n} \sum_{i=1}^n \left(\frac{2}{n} \sum_{j=1}^n \left(9 \exp(2\hat{\beta}_{kA}) + 18 \exp(\hat{\beta}_{kA}) + 9\right) \exp(2 \skew{3.5}\hat{\boldsymbol{\beta}}_{\boldsymbol{kZ}}^T \boldsymbol{Z_j}) \int_{t_q}^{t_s} \frac{dN_{ki}(u)}{\left(S^{(0)}(\skew{3.5}\hat{\boldsymbol{\beta}}_{\boldsymbol{k}}, u)\right)^2} \right. \\[-0.1cm]
&\hphantom{\leq \frac{2}{n} \sum_{i=1}^n \left( i \right.} \left. + \frac{2}{n} \sum_{j=1}^n \left(2\left(\left(\hat{F}_1(t_s \mid A=1, \boldsymbol{Z_j})\right)^2 - \left(\hat{F}_1(t_q \mid A=1, \boldsymbol{Z_j})\right)^2\right) \exp(2 \hat{\beta}_{kA}) \right. \right. \\[-0.3cm]
&\hphantom{\leq \frac{2}{n} \sum_{i=1}^n \left(i + \frac{2}{n} \sum_{j=1}^n \left(i\right.\right.} \left. \left. + 2\left(\left(\hat{F}_1(t_s \mid A=0, \boldsymbol{Z_j})\right)^2 - \left(\hat{F}_1(t_q \mid A=0, \boldsymbol{Z_j})\right)^2\right)\right) \exp(2 \skew{3.5}\hat{\boldsymbol{\beta}}_{\boldsymbol{kZ}}^T \boldsymbol{Z_j}) \int_0^\tau \frac{dN_{ki}(u)}{\left(S^{(0)}(\skew{3.5}\hat{\boldsymbol{\beta}}_{\boldsymbol{k}}, u)\right)^2}\right).
\end{aligned} 
\end{equation}
For the last step, we used again that $\smash{(a+b)^2 \leq 2a^2 + 2b^2}$, as well as the Cauchy-Schwarz inequality and $\smash{(a-b)^2 \leq a^2 - b^2}$ for ${0 \leq b \leq a}$. \\
The second summand in (\ref{eq:lemma4_2}) is further equal to
\begin{align*}
\frac{2}{n} \sum_{i=1}^n \left(\boldsymbol{\hat{\tilde{H}}_{k2}}(t_p) - \boldsymbol{\hat{\tilde{H}}_{k2}}(t_o)\right)^T \! \hat{\mathbf{\Sigma}}_\mathbf{k}^{-1} \tilde{\mathbf{\Sigma}}_\mathbf{ki} \hat{\mathbf{\Sigma}}_\mathbf{k}^{-1} \left(\boldsymbol{\hat{\tilde{H}}_{k2}}(t_p) - \boldsymbol{\hat{\tilde{H}}_{k2}}(t_o)\right),
\end{align*} 
because $\smash{\tilde{\mathbf{\Sigma}}_\mathbf{ki}}$ is symmetric and $\smash{N_{ki}}$ is a one-jump process. Note that
\begin{align*}
\boldsymbol{\hat{\tilde{H}}_{k2}}(t_p) - \boldsymbol{\hat{\tilde{H}}_{k2}}(t_o) &= \frac{1}{n} \int_{t_o}^{t_p} \left(\boldsymbol{\hat{\chi}_{k1, \, A=1}}(u, t_p) - \boldsymbol{\hat{\chi}_{k1, \, A=0}}(u, t_p)\right) \sum_{i=1}^n \frac{dN_{k,i}(u)}{S^{(0)}(\skew{3.5}\hat{\boldsymbol{\beta}}_{\boldsymbol{k}}, u)} \\[-0.1cm]
&\hphantom{= i} - \frac{1}{n} \int_{0}^{t_o} \left(\boldsymbol{\hat{\chi}_{k2, \, A=1}}(u, t_o, t_p) - \boldsymbol{\hat{\chi}_{k2, \, A=0}}(u, t_o, t_p)\right) \sum_{i=1}^n \frac{dN_{k,i}(u)}{S^{(0)}(\skew{3.5}\hat{\boldsymbol{\beta}}_{\boldsymbol{k}}, u)} 
\end{align*}
with 
\begin{gather*}
\begin{multlined}
\boldsymbol{\hat{\chi}_{k1, \, a}}(u, t) = \frac{1}{n} \sum_{i=1}^n \left(\mathbbm{1}\{k=1\} \!\cdot\! \hat{S}(u- \mid a, \boldsymbol{Z_i}) - \hat{F}_1(t \mid a, \boldsymbol{Z_i}) + \hat{F}_1(u \mid a, \boldsymbol{Z_i})\right) \left((a, \boldsymbol{Z}_{\boldsymbol{i}}^T)^T - \boldsymbol{E}(\skew{3.5}\hat{\boldsymbol{\beta}}_{\boldsymbol{k}}, u)\right) \\[-0.3cm]
\cdot \exp(\hat{\beta}_{kA} \!\cdot\! a + \skew{3.5}\hat{\boldsymbol{\beta}}_{\boldsymbol{kZ}}^T \boldsymbol{Z_i}),
\end{multlined} \\
\boldsymbol{\hat{\chi}_{k2, \, a}}(u, s, t) = \frac{1}{n} \sum_{i=1}^n \left(\hat{F}_1(t \mid a, \boldsymbol{Z_i}) - \hat{F}_1(s \mid a, \boldsymbol{Z_i})\right) \left((a, \boldsymbol{Z}_{\boldsymbol{i}}^T)^T - \boldsymbol{E}(\skew{3.5}\hat{\boldsymbol{\beta}}_{\boldsymbol{k}}, u)\right) \exp(\hat{\beta}_{kA} \!\cdot\! a + \skew{3.5}\hat{\boldsymbol{\beta}}_{\boldsymbol{kZ}}^T \boldsymbol{Z_i}).
\end{gather*}
Besides, the product $\smash{\hat{\mathbf{\Sigma}}_\mathbf{k}^{-1} \tilde{\mathbf{\Sigma}}_\mathbf{ki} \hat{\mathbf{\Sigma}}_\mathbf{k}^{-1}}$ is positive definite because of the definitions of $\smash{\hat{\mathbf{\Sigma}}_\mathbf{k}}$ and $\smash{\tilde{\mathbf{\Sigma}}_\mathbf{ki}}$. Since 
\[(\boldsymbol{a}-\boldsymbol{b})^T \mathbf{A} (\boldsymbol{a}-\boldsymbol{b}) \leq 2 \boldsymbol{a}^T \mathbf{A} \boldsymbol{a} + 2 \boldsymbol{b}^T \mathbf{A} \boldsymbol{b} \quad \text{and} \quad \left(\sum_{i_1=1}^n \boldsymbol{a_{i_1}}\right)^T \mathbf{A} \left(\sum_{i_2=1}^n \boldsymbol{a_{i_2}}\right) \leq n \sum_{i=1}^n \boldsymbol{a}_{\boldsymbol{i}}^T \mathbf{A} \boldsymbol{a_i}\]
for a positive (semi-)definite matrix $\mathbf{A}$ and vectors $\boldsymbol{a}$, $\boldsymbol{b}$, $\smash{\boldsymbol{a_1}, \dots, \boldsymbol{a_n}}$, the second summand in (\ref{eq:lemma4_2}) has the upper bound
\begin{align*}
&\frac{2}{n} \sum_{i=1}^n  \left(\frac{2}{n} \int_{t_o}^{t_p} \left(\boldsymbol{\hat{\chi}_{k1, \, A=1}}(u, t_p) - \boldsymbol{\hat{\chi}_{k1, \, A=0}}(u, t_p)\right)^T \hat{\mathbf{\Sigma}}_\mathbf{k}^{-1} \tilde{\mathbf{\Sigma}}_\mathbf{ki} \hat{\mathbf{\Sigma}}_\mathbf{k}^{-1} \left(\boldsymbol{\hat{\chi}_{k1, \, A=1}}(u, t_p) - \boldsymbol{\hat{\chi}_{k1, \, A=0}}(u, t_p)\right) \frac{dN_k(u)}{\left(S^{(0)}(\skew{3.5}\hat{\boldsymbol{\beta}}_{\boldsymbol{k}}, u)\right)^2} \right. \\[-0.1cm]
&\hphantom{\frac{2}{n} \sum_{i=1}^n \left( i \right.} \left. + \frac{2}{n} \int_0^{t_o} \left(\boldsymbol{\hat{\chi}_{k2, \, A=1}}(u, t_o, t_p) - \boldsymbol{\hat{\chi}_{k2, \, A=0}}(u, t_o, t_p)\right)^T \hat{\mathbf{\Sigma}}_\mathbf{k}^{-1} \tilde{\mathbf{\Sigma}}_\mathbf{ki} \hat{\mathbf{\Sigma}}_\mathbf{k}^{-1} \left(\boldsymbol{\hat{\chi}_{k2, \, A=1}}(u, t_o, t_p) - \boldsymbol{\hat{\chi}_{k2, \, A=0}}(u, t_o, t_p)\right) \right. \\[-0.3cm]
&\hphantom{\frac{2}{n} \sum_{i=1}^n \left(i + \frac{2}{n} \int_0^{t_o} i \right.} \left. \cdot \frac{dN_k(u)}{\left(S^{(0)}(\skew{3.5}\hat{\boldsymbol{\beta}}_{\boldsymbol{k}}, u)\right)^2}\right) \\
&\leq \frac{2}{n} \sum_{i=1}^n \left(\frac{4}{n^2} \sum_{j=1}^n \exp(2 \hat{\beta}_{kA} + 2 \skew{3.5}\hat{\boldsymbol{\beta}}_{\boldsymbol{kZ}}^T \boldsymbol{Z_j}) \int_{t_q}^{t_s} \left((1, \boldsymbol{Z}_{\boldsymbol{j}}^T)^T - \boldsymbol{E}(\skew{3.5}\hat{\boldsymbol{\beta}}_{\boldsymbol{k}}, u)\right)^T \hat{\mathbf{\Sigma}}_\mathbf{k}^{-1} \tilde{\mathbf{\Sigma}}_\mathbf{ki} \hat{\mathbf{\Sigma}}_\mathbf{k}^{-1} \left((1, \boldsymbol{Z}_{\boldsymbol{j}}^T)^T - \boldsymbol{E}(\skew{3.5}\hat{\boldsymbol{\beta}}_{\boldsymbol{k}}, u)\right) \frac{dN_k(u)}{\left(S^{(0)}(\skew{3.5}\hat{\boldsymbol{\beta}}_{\boldsymbol{k}}, u)\right)^2} \right. \\[-0.1cm]
&\hphantom{\leq \frac{2}{n} \sum_{i=1}^n \left( i \right.} \left. + \frac{4}{n^2} \sum_{j=1}^n \exp(2 \skew{3.5}\hat{\boldsymbol{\beta}}_{\boldsymbol{kZ}}^T \boldsymbol{Z_j}) \int_{t_q}^{t_s} \left((0, \boldsymbol{Z}_{\boldsymbol{j}}^T)^T - \boldsymbol{E}(\skew{3.5}\hat{\boldsymbol{\beta}}_{\boldsymbol{k}}, u)\right)^T \hat{\mathbf{\Sigma}}_\mathbf{k}^{-1} \tilde{\mathbf{\Sigma}}_\mathbf{ki}\hat{\mathbf{\Sigma}}_\mathbf{k}^{-1} \left((0, \boldsymbol{Z}_{\boldsymbol{j}}^T)^T - \boldsymbol{E}(\skew{3.5}\hat{\boldsymbol{\beta}}_{\boldsymbol{k}}, u)\right) \frac{dN_k(u)}{\left(S^{(0)}(\skew{3.5}\hat{\boldsymbol{\beta}}_{\boldsymbol{k}}, u)\right)^2} \right. \\[-0.1cm]
&\hphantom{\leq \frac{2}{n} \sum_{i=1}^n \left( i \right.} \left. + \frac{4}{n^2} \sum_{j=1}^n \left(\left(\hat{F}_1(t_s \mid A=1, \boldsymbol{Z_j})\right)^2 - \left(\hat{F}_1(t_q \mid A=1, \boldsymbol{Z_j})\right)^2\right) \exp(2 \hat{\beta}_{kA} + 2 \skew{3.5}\hat{\boldsymbol{\beta}}_{\boldsymbol{kZ}}^T \boldsymbol{Z_j}) \right. \\[-0.2cm]
&\hphantom{\leq \frac{2}{n} \sum_{i=1}^n \left( i + \frac{4}{n^2} \sum_{j=1}^n \left( i\right. \right.} \left. \cdot \int_0^\tau \left((1, \boldsymbol{Z}_{\boldsymbol{j}}^T)^T - \boldsymbol{E}(\skew{3.5}\hat{\boldsymbol{\beta}}_{\boldsymbol{k}}, u)\right)^T \hat{\mathbf{\Sigma}}_\mathbf{k}^{-1} \tilde{\mathbf{\Sigma}}_\mathbf{ki} \hat{\mathbf{\Sigma}}_\mathbf{k}^{-1} \left((1, \boldsymbol{Z}_{\boldsymbol{j}}^T)^T - \boldsymbol{E}(\skew{3.5}\hat{\boldsymbol{\beta}}_{\boldsymbol{k}}, u)\right) \frac{dN_k(u)}{\left(S^{(0)}(\skew{3.5}\hat{\boldsymbol{\beta}}_{\boldsymbol{k}}, u)\right)^2} \right. \\[-0.1cm]
&\hphantom{\leq \frac{2}{n} \sum_{i=1}^n \left( i \right.} \left. + \frac{4}{n^2} \sum_{j=1}^n \left(\left(\hat{F}_1(t_s \mid A=0, \boldsymbol{Z_j})\right)^2 - \left(\hat{F}_1(t_q \mid A=0, \boldsymbol{Z_j})\right)^2\right) \exp(2 \skew{3.5}\hat{\boldsymbol{\beta}}_{\boldsymbol{kZ}}^T \boldsymbol{Z_j}) \right. \\[-0.5cm]
&\hphantom{\leq \frac{2}{n} \sum_{i=1}^n \left( i + \frac{4}{n^2} \sum_{j=1}^n \left( i \right. \right.} \left. \cdot \int_0^\tau \left((0, \boldsymbol{Z}_{\boldsymbol{j}}^T)^T - \boldsymbol{E}(\skew{3.5}\hat{\boldsymbol{\beta}}_{\boldsymbol{k}}, u)\right)^T \hat{\mathbf{\Sigma}}_\mathbf{k}^{-1} \tilde{\mathbf{\Sigma}}_\mathbf{ki} \hat{\mathbf{\Sigma}}_\mathbf{k}^{-1} \left((0, \boldsymbol{Z}_{\boldsymbol{j}}^T)^T - \boldsymbol{E}(\skew{3.5}\hat{\boldsymbol{\beta}}_{\boldsymbol{k}}, u)\right) \frac{dN_k(u)}{\left(S^{(0)}(\skew{3.5}\hat{\boldsymbol{\beta}}_{\boldsymbol{k}}, u)\right)^2}\right).
\end{align*}
For the last two terms, we used that for ${0 \leq b \leq a}$, $\smash{(a-b)^2 \leq a^2 - b^2}$. Note also that it is possible to extend the integral limits here because of the positive definiteness of $\smash{\hat{\mathbf{\Sigma}}_\mathbf{k}^{-1} \tilde{\mathbf{\Sigma}}_\mathbf{ki} \hat{\mathbf{\Sigma}}_\mathbf{k}^{-1}}$. \\
Together with (\ref{eq:lemma4_3}), it follows finally that
\[\sum_{i=1}^n \left(X_{n,i}^{(k)}(t_p) - X_{n,i}^{(k)}(t_o)\right)^2 \leq 36 \left(L_n^{(k)}(t_s) - L_n^{(k)}(t_1)\right) \]
${\forall k \in \{1, \dots, K\}}$.
\end{proof}

\bibliography{bib_Theory}{}
\bibliographystyle{apalike}

\end{document}